\documentclass[a4paper,12pt]{amsart}
\usepackage{amssymb}
\usepackage{ifthen}
\usepackage{graphicx}
\usepackage{multicol}        
\usepackage{epsfig}

\usepackage{latexsym}
\usepackage{enumerate}

\nonstopmode \numberwithin{equation}{section}
\newtheorem{thm}[equation]{Theorem}
\newtheorem{cor}[equation]{Corollary}
\newtheorem{lem}[equation]{Lemma}
\newtheorem{prop}[equation]{Proposition}

\newtheorem{conj}{Conjecture}

\theoremstyle{definition}

\newtheorem{prob}[equation]{Problem}

\newenvironment{rem}{%
\bigskip
\noindent \textsl{{\sl Remark. }}}{\bigskip}
\newenvironment{rems}{%
\bigskip
\noindent \textsl{{\sl Remarks. }}}{\bigskip}

\newcounter {own}
\def\theown {\thesection       .\arabic{own}}

\newenvironment{pf}[1][]{%
 \vskip 3mm
 \noindent
 \ifthenelse{\equal{#1}{}}%
  {{\slshape Proof. }}%
  {{\slshape #1.} }%
 }%
{\qed\bigskip}

\newcounter{alphabet}
\newcounter{tmp}

\newcommand{\ID}{{\mathbb D}}

\newcommand{\IC}{{\mathbb C}}

\def\be{\begin{equation}}
\def\ee{\end{equation}}

\newcommand{\bee}{\begin{enumerate}}
\newcommand{\eee}{\end{enumerate}}

\newcommand{\blem}{\begin{lem}}
\newcommand{\elem}{\end{lem}}
\newcommand{\bthm}{\begin{thm}}
\newcommand{\ethm}{\end{thm}}
\newcommand{\bcor}{\begin{cor}}
\newcommand{\ecor}{\end{cor}}
\newcommand{\beg}{\begin{examp}}
\newcommand{\eeg}{\end{examp}}
\newcommand{\begs}{\begin{examples}}
\newcommand{\eegs}{\end{examples}}
\newcommand{\bdefe}{\begin{defin}}
\newcommand{\edefe}{\end{defin}}
\newcommand{\bprob}{\begin{prob}}
\newcommand{\eprob}{\end{prob}}
\newcommand{\bei}{\begin{itemize}}
\newcommand{\eei}{\end{itemize}}

\newcommand{\bcon}{\begin{conj}}
\newcommand{\econ}{\end{conj}}
\newcommand{\bcons}{\begin{conjs}}
\newcommand{\econs}{\end{conjs}}
\newcommand{\bprop}{\begin{prop}}
\newcommand{\eprop}{\end{prop}}
\newcommand{\br}{\begin{rem}}
\newcommand{\er}{\end{rem}}
\newcommand{\brs}{\begin{rems}}
\newcommand{\ers}{\end{rems}}
\newcommand{\bo}{\begin{obser}}
\newcommand{\eo}{\end{obser}}
\newcommand{\bos}{\begin{obsers}}
\newcommand{\eos}{\end{obsers}}
\newcommand{\bpf}{\begin{pf}}
\newcommand{\epf}{\end{pf}}
\newcommand{\ba}{\begin{array}}
\newcommand{\ea}{\end{array}}
\newcommand{\beq}{\begin{eqnarray}}
\newcommand{\beqq}{\begin{eqnarray*}}
\newcommand{\eeq}{\end{eqnarray}}
\newcommand{\eeqq}{\end{eqnarray*}}

\newcommand{\ds}{\displaystyle}

\newcounter{minutes}
\divide\time by 60
\newcounter{hours}
\multiply\time by 60 \addtocounter{minutes}{-\time}

\begin{document}
\bibliographystyle{amsplain}
\title{Region of variability for  functions with positive real part}
\author{S. Ponnusamy}
\address{S. Ponnusamy, Department of Mathematics,
Indian Institute of Technology Madras, Chennai-600 036, India.}
\email{samy@iitm.ac.in}
\author{A. Vasudevarao}
\address{A. Vasudevarao, Department of Mathematics,
Indian Institute of Technology Madras, Chennai-600 036, India.}
\email{alluvasudevarao@gmail.com}


\date{
\texttt{File:~\jobname .tex,
          printed: \number\year-\number\month-\number\day,
          \thehours.\ifnum\theminutes<10{0}\fi\theminutes}
}
\begin{abstract}

For  $\gamma\in\IC$ such that $|\gamma|<\pi/2$ and $0\leq\beta<1$, let
${\mathcal P}_{\gamma,\beta} $ denote
the class of all analytic functions $P$ in the unit disk $\mathbb{D}$
with $P(0)=1$ and
$$ {\rm Re\,} \left (e^{i\gamma}P(z)\right)>\beta\cos\gamma \quad \mbox{ in  ${\mathbb D}$}.
$$
For any fixed $z_0\in\mathbb{D}$ and
$\lambda\in\overline{\mathbb{D}}$,
we shall determine the region of variability $V_{\mathcal{P}}(z_0,\lambda)$ for
$\int_0^{z_0}P(\zeta)\,d\zeta$ when $P$ ranges over the class
$$
\mathcal{P}(\lambda)   =
\left\{ P\in{\mathcal P}_{\gamma,\beta} :\,
P'(0)=2(1-\beta)\lambda e^{-i\gamma}\cos\gamma
\right\}.
$$
As a consequence, we present the region of variability for some subclasses
of univalent functions.  We also graphically illustrate the
region of variability for several sets of parameters.

\end{abstract}

\thanks{}

\maketitle
\pagestyle{myheadings}
\markboth{S. Ponnusamy and  A. Vasudevarao }{Regions of variability}

\section{Introduction}\label{sec-01}

We denote by $\mathcal{H}$ the class of analytic functions in the unit disk
$\ID=\{z\in\IC:\,|z|<1\}$, and think of
$\mathcal{H}$ as a topological vector space endowed with the
topology of uniform convergence over compact subsets of
$\mathbb{D}$. We consider  the subclass of functions
$\phi\in\mathcal{H}$ with $\phi(0)=0=\phi'(0)-1$ such that $\phi$ maps
$\mathbb{D}$ univalently onto a domain
that is starlike (with respect to the origin). That is,
$t\phi(z)\in\phi(\mathbb{D})$ for each  $t\in [0,1]$. We denote the class
of such functions by $\mathcal{S}^*$. Analytically,  each $\phi\in {\mathcal S}^*$ is characterized
by the condition
$${\rm Re\,}\left (\frac{z\phi '(z)}{\phi (z)}\right )>0, \quad z\in \ID.
$$
Functions in ${\mathcal S}^*$ are referred to as {\em starlike functions}.
A function $\phi\in\mathcal{H}$ with
$\phi(0)=0=\phi'(0)-1$ is said to belong to $\mathcal{C}$ if and only if
$\phi(\mathbb{D})$ is a convex domain. It is well-known that $\phi\in\mathcal{C}$
if and only if $z\phi'\in{\mathcal S}^*$. Functions in ${\mathcal C}$ are
referred to as convex functions.

Let ${\mathcal P}_{\gamma,\beta}$ denote the class of functions
$P\in\mathcal{H}$ with $P(0)=1$ and
$${\rm Re\,} \left (e^{i\gamma}P(z)\right)>\beta\cos\gamma \quad \mbox{ in  ${\mathbb D}$},
$$
for some $\beta$ with $\beta<1$ and $\gamma\in\mathbb{C}$ with $|\gamma|<\pi/2$.
Let $\mathcal{A}$ denote the  class of functions $f$ in $\mathcal{H}$ such that
$f(0)=0=f'(0)-1$.
When $P(z)=\frac{zf'(z)}{f(z)}$ and $\beta=0$, the class $\mathcal{P}_{\gamma,\beta}$
becomes
$$
\mathcal{S}^{\gamma}(0)=\left\{f\in\mathcal{A}\colon
 {\rm Re\,}\left(e^{i\gamma}\frac{zf'(z)}{f(z)}\right)>0\quad\mbox{ in }\mathbb{D}\right\}
$$
for some $\gamma$ with $|\gamma|<\pi/2$. Functions in $\mathcal{S}^{\gamma}(0)$
are known to be univalent in $\mathbb{D}$ and $\mathcal{S}^{0}(0)\equiv\mathcal{S}^*$.
Functions in $\mathcal{S}^{\gamma}(0)$ are called spirallike functions (see \cite{Spacek-33}).

\section{Preliminary Investigation about the class
$\mathcal{P}_{\gamma,\beta}$}\label{sec-02}

Herglotz representation for analytic functions with positive real part in
$\mathbb{D}$ shows that if $P\in\mathcal{P}_{\gamma,\beta}$, then there exists a
unique positive unit measure $\mu$ on $(-\pi,\pi]$ such that
$$P(z)=
\int_{-\pi}^{\pi}\frac{1+[1-2\beta e^{-i\gamma}\cos\gamma]ze^{-it}}{1-ze^{-it}}\,d\mu(t).
$$
Let ${\mathcal B}_0$ be the class of analytic functions $\omega$
in $\mathbb{D}$ such that $| \omega(z)|<1$ in $\mathbb{D}$ and
$\omega(0)=0$. Then it is a simple exercise to see that for each
$P\in {\mathcal P}_{\gamma,\beta}$ there exists an $\omega_P \in {\mathcal B}_0$
such that
\be\label{pvdev-10-eq2}
\omega_P(z)=\frac{e^{i\gamma}P(z)-e^{i\gamma}}
{e^{i\gamma}P(z)-(2\beta\cos\gamma- e^{-i\gamma})},\quad z\in\mathbb{D},
\ee
and conversely. Clearly, we have
$$
P'(0)=2e^{-i\gamma}\omega'_P(0)(1-\beta)\cos\gamma.
$$
Suppose that $P\in {\mathcal P}_{\gamma,\beta}$.
Then, because $|\omega'_P(0)|\leq 1$,  by  the classical Schwarz lemma
(see for example \cite{Dinnen-book,Du,samy-book3,samy-herb}) we may let
$$P'(0)= 2\lambda e^{-i\gamma}(1-\beta)\cos\gamma
$$
for some $\lambda\in\overline{\mathbb{D}}$, with $\omega'_P(0)=\lambda$.
Using (\ref{pvdev-10-eq2}), one can compute
\be\label{pvdev-10-eq2b}
\frac{\omega''_P(0)}{2}=\frac{e^{i\gamma}P''(0)}{4(1-\beta)\cos\gamma}-{\lambda}^2.
\ee
Also if we let
$$g(z)=\left\{\ba{ll}
\ds\frac{\frac{\omega_P(z)}{z}-\lambda}
{1-\overline{\lambda}\frac{\omega_P(z)}{z}} & \mbox{ for } |\lambda|< 1\\ [6mm]
\ds 0 & \mbox{  for }  |\lambda|= 1,
\ea \right.
$$
then we see that
$$g'(0)=\left\{\ba{ll}
\ds\left.\frac{1}{1-|\lambda|^2}\left(\frac{\omega_P(z)}{z}\right)'\right |_{z=0}
=\frac{1}{1-|\lambda|^2}\frac{\omega''_P(0)}{2}
& \mbox {  for  } |\lambda|< 1 \\ [6mm]
\ds 0 & \mbox{  for  } |\lambda|= 1.
\ea\right.
$$
By the Schwarz lemma, $|g(z)|\leq |z|$ and $|g'(0)|\leq 1$.
Equality holds in both the cases
if and only if $g(z)=e^{i\alpha}z$ for some $\alpha\in\mathbb{R}$.
The condition  $|g'(0)|\leq 1$ shows that there exists
an $a\in\overline{\mathbb{D}}$ such that $g'(0)=a$.
In view of (\ref{pvdev-10-eq2b}) we may represent $P''(0)$ as
\be\label{pvdev-10-eq25}
P''(0)=4(1-\beta)[(1-|\lambda|^2)a+ {\lambda}^2]e^{-i\gamma}\cos\gamma
\ee
for some $a\in\overline{\mathbb{D}}$.
Consequently, for
$\lambda\in\overline{\mathbb{D}}=\{z\in\mathbb{C}:\,|z|\leq 1\}$  and
$z_0\in\mathbb{D}$ fixed, it is natural to introduce (for convenience with the
notation $\mathcal{P}(\lambda)$
instead of ${\mathcal P}_{\gamma,\beta}(\lambda)$)
\beqq
\mathcal{P}(\lambda)  & = &
\left\{\frac{}{} P\in{\mathcal P}_{\gamma,\beta} :\,
P'(0)=2(1-\beta)e^{-i\gamma}\lambda\cos\gamma
\right\}\\
V_{\mathcal{P}}(z_0,\lambda) &  = &  \left\{\int_0^{z_0}P(\zeta)\,d\zeta:
\, P\in {\mathcal P}(\lambda)\right\}.
\eeqq
Obviously, each $f \in \mathcal{P}(\lambda) $ has to satisfy the condition
(\ref{pvdev-10-eq25}) for some  $a\in\overline{\mathbb{D}}$ and so we do not need to include
it in the defintion  of $\mathcal{P}(\lambda) $.

For each fixed $z_0\in \ID$,
using extreme function theory, it has been shown by
Grunsky \cite[Theorem~10.6]{Du} that the region of variability of
$$V(z_{0})=\left \{\log\frac{f(z_{0})}{z_{0}}:\ f\in\mathcal{S} \right \}
$$
is precisely a closed disk, where $\mathcal{S} =\{f\in \mathcal{A}:\,
\mbox{$f$ is univalent in $\ID$} \}.$ It is also well-known that the region of variability
$$V(z_0)=\{\log \phi ' (z_0):\,\phi \in {\mathcal C}\}
$$
is the set $\{\log (1-z)^{-2}:\,|z|\leq |z_0|\}$.
Recently,  several authors have studied region of variability problems
for various subclasses of univalent functions in  $\mathcal{H}$, see
\cite{Paatero,Pinchuk,samy-vasudev2,samy-vasudev3,
samy-vasudev-yan3,samy-vasudev-vuorinen1,samy-vasudev-vuorinen2,Yanagihara1,Yanagihara2}.

The main aim of this paper is to determine
the region of variability of $V_{\mathcal{P}}(z_0,\lambda)$ for $\int_0^{z_0}P(\zeta)\,d\zeta$ when
$P$ ranges over the class ${\mathcal P}(\lambda)$.
In Section \ref{sec-03}, we present some basic properties of $V_{\mathcal{P}}(z_0,\lambda)$
whereas in  Section \ref{sec-04}, we investigate the growth condition for functions in ${\mathcal P}(\lambda)$.
The precise geometric description of the set $V_{\mathcal{P}}(z_0,\lambda)$ is established
in Theorem \ref{pvdev-10-th1} in  Section \ref{sec-05}.
Two interesting special cases are presented in  Section \ref{sec-06}.
Finally, in Section \ref{sec-07}, we graphically illustrate the
region of variability for several sets of parameters.


\section{Basic properties of $V_{\mathcal{P}}(z_0,\lambda)$}\label{sec-03}

For a positive integer $p$, let
$({\mathcal S}^*)^p=\{f=f_0^p:\, f_0\in {\mathcal S}^* \}.
$
A sufficient condition (see \cite{Yanagihara1}) for an analytic function $f$
in ${\mathbb D}$ with $f(z) = z^p + \cdots $,  to be in $({\mathcal S}^*)^p$
is that $f$ satisfies
\be\label{pvdev-10-eq2c}
{\rm Re} \,  \left( 1+ z \frac{f''(z)}{f'(z)} \right)> 0 , \quad z \in {\mathbb D}.
\ee
This fact will be used in the following result.

\bprop\label{pvdev-10-pro01}
We have
\bee
\item[\rm{(1)}] $V_{\mathcal{P}}(z_0,\lambda)$ is a compact subset of $\mathbb{C}$.
\item[\rm{(2)}] $V_{\mathcal{P}}(z_0,\lambda)$ is a convex subset of $\mathbb{C}$.
\item[\rm{(3)}] For $|\lambda|=1$ or $z_0=0$,
\be\label{pvdev-10-eq20}
V_{\mathcal{P}}(z_0,\lambda)=\left \{z_0-2(1-\beta)e^{-i\gamma}\cos\gamma
\left(z_0+\frac{1}{\lambda}\log(1-\lambda z_0)\right)\right \}.
\ee
\item[\rm{(4)}] For $|\lambda|<1$ and $z_0\in\mathbb{D}\setminus\{0\}$, $V_{\mathcal{P}}(z_0,\lambda)$ has
$$z_0-2(1-\beta)e^{-i\gamma}\cos\gamma
\left(z_0+\frac{1}{\lambda}\log(1-\lambda z_0)\right)
$$
as an interior point.
\eee
\eprop

\bpf
(1) Since ${\mathcal P}(\lambda)$ is a
compact subset of $\mathcal{H}$, it follows that
$V_{\mathcal{P}}(z_0,\lambda)$ is also compact.

(2) If $p_1,p_2\in {\mathcal P}(\lambda)$ and $0\leq t \leq 1$, then
the function
$$P_t(z)=(1-t)p_1(z)+tp_2(z)
$$
is evidently in ${\mathcal P}(\lambda)$. Also, because of the
representation of $P_t$, we see easily that  the set $V_{\mathcal{P}}(z_0,\lambda)$ is
convex.

(3) If $z_0=0$, (\ref{pvdev-10-eq20}) trivially holds. If $|\lambda|=1$,
then from our earlier observation $\omega_P(z)=\lambda z$  and so,
$P\in{\mathcal P}(\lambda)$ defined by (\ref{pvdev-10-eq2}) takes the form
$$P(z)=\frac{1+\lambda z \left[2(1-\beta)e^{-i\gamma}\cos\gamma-1\right]}{1-\lambda z}
$$
or equivalently,
$$P(z)=1-2(1-\beta)e^{-i\gamma}\cos\gamma\left(1-\frac{1}{1-\lambda z}\right).
$$
Consequently,
$$V_{\mathcal{P}}(z_0,\lambda)
=\left \{z_0-2(1-\beta)e^{-i\gamma}\cos\gamma
\left(z_0+\frac{1}{\lambda}\log(1-\lambda z_0)\right)\right \}.
$$

(4) For $|\lambda|<1$ and $a\in\overline{\mathbb{D}}$, we let
\be\label{pvdev-10-delta}
\delta(z,\lambda) = \frac{z+\lambda}{1+\overline{\lambda}z},
\ee
and in order to get the extremal function in ${\mathcal P}(\lambda)$,
we define
\be\label{pvdev-10-extremal}
H_{a,\lambda}(z) = 1+2(1-\beta)e^{-i\gamma}\cos\gamma
\frac{\delta(az, \lambda)z}{1-\delta(az, \lambda)z}.
\ee
Clearly $H_{a,\lambda}(0)=1$.
Since $\delta(az,\lambda)$ lies in the unit disk $\mathbb{D}$ and $\varphi(w)=w/(1-w)$
maps $|w|<1$ onto ${\rm Re\,}\varphi(w)>-1/2$, we obtain that
$${\rm Re\,}\left(e^{i\gamma}H_{a,\lambda}(z)\right)>\beta\cos\gamma ~\mbox { in }\mathbb{D}.
$$
Also, from (\ref{pvdev-10-extremal}), we have the normalization condition
$$ H'_{a,\lambda}(0)  =  2(1-\beta)e^{-i\gamma}\lambda\cos\gamma .
$$
Thus, $H_{a,\lambda}\in {\mathcal P}(\lambda)$. We observe that
\be\label{pvdev-10-eq5}
\omega_{H_{a,\lambda}}(z)=z\delta(az,\lambda).
\ee

We claim that the mapping
$${\mathbb D} \ni a\mapsto \int_0^{z_0}H_{a,\lambda}(\zeta)\,d\zeta
$$
is a non-constant analytic function of $a$
for each fixed $z_0 \in {\mathbb D} \backslash \{ 0 \}$ and
$\lambda\in\mathbb{D}$. To see this, we introduce
\begin{eqnarray*}
h(z) & = & \left .\frac{3e^{i\gamma}}{2(1-\beta)(1-|\lambda|^2)\cos\gamma}
 \frac{\partial}{\partial a}\left\{\frac{}{}\int_0^{z}H_{a,\lambda}
(\zeta)\,d\zeta\right\}\right |_{a=0}
\end{eqnarray*}
so that
\begin{eqnarray*}
h(z) & = & \left .\frac{3}{1-|\lambda|^2}
\frac{\partial}{\partial a}\left\{\frac{}{}\int_0^{z}
\frac{\delta(a\zeta,\lambda)\zeta}{1-\delta(a\zeta,\lambda)\zeta}
\,d\zeta\right\}\right |_{a=0}.
\end{eqnarray*}
A computation gives
\begin{eqnarray*}
h(z) & = & \left .
3\frac{\partial}{\partial a}\left\{\frac{}{}\int_0^{z}
\frac{{\zeta}^2}{(1-\lambda \zeta)^2}
\frac{d\zeta}{(1-a\delta(a\zeta,\lambda)\zeta)^2}
\right\}\right |_{a=0}
\end{eqnarray*}
which clearly shows that
\begin{eqnarray*}
h(z) & = &  3\int_0^z \frac{{\zeta}^2}{(1-\lambda\zeta)^2}\,d\zeta=z^3+\cdots
\end{eqnarray*}
from which it is easy to see that
$${\rm Re} \,\left\{\frac{zh''(z)}{h'(z)}\right\}=
2{\rm Re} \,\left\{\frac{1}{1-\lambda z}\right\}>\frac{2}{1+|\lambda|}\geq 1,
\quad z\in\mathbb{D}.
$$
By (\ref{pvdev-10-eq2c}), there exists a function $h_0\in
{\mathcal S}^*$  with $h=h_0^3$. The univalence of $h_0$ together with
the condition $h_0(0)=0$ implies that $h(z_0)\neq 0$ for $z_0 \in {\mathbb
D}\setminus \{0\}$. Consequently, the mapping ${\mathbb D} \ni
a\mapsto \int_0^{z_0} H_{a,\lambda}(\zeta)\,d\zeta$ is a non-constant analytic
function of $a$ and hence, it is an open mapping. Thus,
$V_{\mathcal{P}}(z_0,\lambda)$ contains the open set
$$\left\{\int_0^{z_0} H_{a,\lambda}(\zeta)\,d\zeta:\, |a|<1\right\}.
$$
In particular,
$$\int_0^{z_0} H_{0,\lambda}(\zeta)\,d\zeta=z_0-2(1-\beta)e^{-i\gamma}\cos\gamma
\left(z_0+\frac{1}{\lambda}\log(1-\lambda z_0)\right)
$$
is an interior point of

\vspace{6pt}

\hfill$\displaystyle \left\{\int_0^{z_0}H_{a,\lambda}(\zeta)\,d\zeta:\,a\in\mathbb{D}\right\}
\subset V_{\mathcal{P}}(z_0,\lambda).
$\hfill
\epf

We remark that, since $V_{\mathcal{P}}(z_0,\lambda)$ is a compact convex subset of
$\mathbb{C}$ and has nonempty interior, the boundary
$\partial{V_{\mathcal{P}}(z_0,\lambda)}$ is a Jordan curve and $V_{\mathcal{P}}(z_0,\lambda)$
is the union of $\partial{V_{\mathcal{P}}(z_0,\lambda)}$ and its inner domain.


\section{Growth condition for functions in $\mathcal{P}(\lambda)$} \label{sec-04}

\bprop\label{pvdev-10-pro1} For $P\in {\mathcal P}(\lambda)$  with
$\lambda\in\mathbb{D}$, we have
\be\label{pvdev-10-eq8}
\left|P(z)-c(z,\lambda)\right|\leq r(z,\lambda), \quad z\in\mathbb{D},
\ee
where
\begin{eqnarray*}
c(z,\lambda)& = & \frac{(1+\lambda z(e^{-i\gamma}-2\beta\cos\gamma)e^{-i\gamma})
(1-\overline{\lambda}\overline{z})}{(1-|z|^2)(1+|z|^2-2{\rm Re\,}({\lambda}z))}\\
& & \quad +\frac{|z|^2(\overline{z}-\lambda)
\left(\overline{\lambda}+z(e^{-i\gamma}-2\beta\cos\gamma)e^{-i\gamma}\right)}
{(1-|z|^2)(1+|z|^2-2{\rm Re\,}({\lambda} z))}, \\
r(z,\lambda) & = & \frac{2(1-|\lambda|^2)(1-\beta)|z|^2\cos\gamma}
{(1-|z|^2)(1+|z|^2-2{\rm Re\,}({\lambda} z))}.
\end{eqnarray*}
For each $z\in\mathbb{D}\setminus\{0\}$, equality holds if and
only if $P=H_{e^{i\theta},\lambda}$ for some
$\theta\in\mathbb{R}$.
\eprop \bpf Let $P\in{\mathcal P}(\lambda)$. Then there exists
$\omega_P \in {\mathcal B}_0$ satisfying
(\ref{pvdev-10-eq2}). As observed in Section \ref{sec-02} $(|g(z)|\leq |z|)$, we have
\be\label{pvdev-10-eq9}
\left|\frac{\frac{\omega_P(z)}{z}-\lambda}
{1-\overline{\lambda}\frac{\omega_P(z)}{z}}\right|\leq|z|, \quad z\in\mathbb{D}.
\ee
From   (\ref{pvdev-10-eq2}) this is
equivalent to
\be\label{pvdev-10-eq10}
\left|\frac{P(z)-A(z,\lambda)}{P(z)+B(z,\lambda)}
\right| \leq |z|\, |\tau(z,\lambda)|,
\ee
where
\be\label{pvdev-10-eq11}
\left\{
\ba{lll}
A(z,\lambda)& =& \ds \frac{1+e^{-i\gamma}\lambda z(e^{-i\gamma}-2\beta\cos\gamma)}
{1-\lambda z}\\
B(z,\lambda) & = & \ds \frac{\overline{\lambda}+
e^{-i\gamma}z(e^{-i\gamma}-2\beta\cos\gamma)}
{z-\overline{\lambda}}\\
\tau(z,\lambda)& = & \ds \frac{z-\overline{\lambda}}{1-\lambda z}.
\ea
\right.
\ee
A simple calculation shows that the inequality
(\ref{pvdev-10-eq10}) is equivalent to
\be\label{pvdev-10-eq12}
\left|P(z)-\frac{A(z,\lambda)+|z|^2\,|\tau(z,\lambda)|^2 B(z,\lambda)}
{1-|z|^2\,|\tau(z,\lambda)|^2}\right| \leq
\frac{|z|\,|\tau(z,\lambda)|\,|A(z,\lambda)+B(z,\lambda)|\,}{1-|z|^2\,|\tau(z,\lambda)|^2}.
\ee
Using (\ref{pvdev-10-eq11}) we can easily see that
\begin{eqnarray}
\label{pvdev-10-eq13} 1-|z|^2\, |\tau(z,\lambda)|^2 & = &
\frac{(1-|z|^2)(1+|z|^2-2{\rm Re\,}(\lambda z))}{|1-\lambda
z|^2}\nonumber,
\\
\label{pvdev-10-eq14}
A(z,\lambda)+ B(z,\lambda)
& = & \frac{2(1-|\lambda|^2)(1-\beta)(\cos\gamma)e^{-i\gamma}z}{(1-\lambda z)
(z-\overline{\lambda})}\nonumber
\end{eqnarray}
and
\begin{eqnarray*}
A(z,\lambda) &+&|z|^2|\tau(z,\lambda)|^2 B(z,\lambda)\\
&=&  \frac{(1+\lambda z(e^{-i\gamma}-2\beta\cos\gamma)e^{-i\gamma})
(1-\overline{\lambda}\overline{z})}{|1-\lambda z|^2}\\
& & \quad +\frac{|z|^2(\overline{z}-\lambda)\left(\overline{\lambda}+z(e^{-i\gamma}
-2\beta\cos\gamma)e^{-i\gamma}\right)}{|1-\lambda z|^2}.
\end{eqnarray*}
Thus, by a simple computation, we see that
$$\frac{A(z,\lambda)+|z|^2|\tau(z,\lambda)|^2 B(z,\lambda)}
{1-|z|^2|\tau(z,\lambda)|^2} = c(z,\lambda)
$$
and
$$\frac{|z|\,|\tau(z,\lambda)|\,|A(z,\lambda)+B(z,\lambda)|}{1-|z|^2|\tau(z,\lambda)|^2}
=r(z,\lambda).
$$
Now the inequality (\ref {pvdev-10-eq8}) follows from these
equalities and (\ref{pvdev-10-eq12}).

It is easy to see that the equality occurs  in
(\ref {pvdev-10-eq8}) for a $z\in\mathbb{D}$, when $P=H_{e^{i\theta},\lambda}$ for some
$\theta\in\mathbb{R}$. Conversely if the equality occurs for some
$z\in\mathbb{D}\setminus\{0\}$ in (\ref {pvdev-10-eq8}), then the
equality must hold in (\ref{pvdev-10-eq9}). Thus from the Schwarz
lemma there exists a $\theta\in\mathbb{R}$ such that
$\omega_P(z)=z\delta(e^{i\theta}z,\lambda)$ for all $z\in\mathbb{D}$. This
implies $P=H_{e^{i\theta},\lambda}$.
\epf

The choice of $\lambda=0$ gives the following result which may deserve
a special mention.

\bcor\label{pvdev-10-cor01}
For $P\in{\mathcal P}(0)$ we have
\be\label{pvdev-10-eq24}
\left|P(z)-\frac{1+(1-2\beta)|z|^4)}{1-|z|^4}\right|\leq
\frac{2(1-\beta)|z|^2}{1-|z|^4}, \quad z\in\mathbb{D}.
\ee

For each
$z\in\mathbb{D}\setminus\{0\}$, equality holds if and only if
$P=H_{e^{i\theta},0}$ for some $\theta\in\mathbb{R}$.
\ecor

\bcor\label{pvdev-10-cor1}
Let $\gamma:\,z(t)$, $0\leq t\leq 1$, be
a $C^1$-curve in $\mathbb{D}$ with $z(0)=0$ and $z(1)=z_0$. Then
we have
$$V_{\mathcal{P}}(z_0,\lambda)\subset
\{w\in\mathbb{C}:\,|w-C(\lambda, \gamma) |\leq R(\lambda, \gamma)\},
$$
where
$$C(\lambda, \gamma)=\int_0^1 c(z(t),\lambda)z'(t)\,dt ~\mbox{ and }~ R(\lambda,
\gamma)=\int_0^1 r(z(t),\lambda)|z'(t)|\,dt.
$$
\ecor \bpf
Proof follows as in \cite{samy-vasudev-vuorinen2}.
\epf


For the proof of our next result, we need the following lemma.

\blem\label{pvdev-10-lem01}
For $\theta\in\mathbb{R}$ and $\lambda\in\mathbb{D}$, the function
$$G(z)=\int_0^z \frac{ e^{i\theta}{\zeta}^2 }
{\{1+(\overline{\lambda}e^{i\theta}-\lambda)\zeta-e^{i\theta}{\zeta}^2\}^2}\,
d\zeta, \quad z\in\mathbb{D},
$$
has zero of order three at the origin and no zeros elsewhere in
$\mathbb{D}$. Furthermore there exists a starlike univalent
function $G_0$ in $\mathbb{D}$ such that
$G=(1/3)e^{i\theta}G^3_0$ and $G_0(0)= G'_0(0)-1=0$.
\elem \bpf
For a proof, we refer to \cite[Lemma 3.4]{samy-vasudev-vuorinen2} with $\beta=1$ there.
\epf

\bprop
Let $z_0\in\mathbb{D}\setminus \{0\}$. Then for $\theta\in(-\pi,\pi]$ we have
$$\int_0^{z_0}H_{e^{i\theta},\lambda}(\zeta)\,d\zeta\in\partial V_{\mathcal{P}}(z_0,\lambda).
$$
Furthermore if $\int_0^{z_0}P(\zeta)\,d\zeta= \int_0^{z_0}H_{e^{i\theta},\lambda}
(\zeta)\,d\zeta$ for some $P\in {\mathcal P}(\lambda)$ and $\theta\in(-\pi,\pi]$,
then $P= H_{e^{i\theta},\lambda}$.
\label{pvdev-10-pro2}
\eprop\bpf
From (\ref{pvdev-10-extremal}) we have
\begin{eqnarray*}
H_{a, \lambda}(z)
& = & \frac{1+[2(1-\beta)(\cos\gamma)e^{-i\gamma}-1]\delta(az, \lambda)z}
{1-\delta(az, \lambda)z} \\
& = & \frac{1+\overline{\lambda}az+(\lambda z+az^2)(2(1-\beta)(\cos\gamma)e^{-i\gamma}-1)}
{1+(\overline{\lambda}a-\lambda)z-az^2}.
\end{eqnarray*}
Using (\ref{pvdev-10-eq11}) we compute
$$H_{a, \lambda}(z)-A(z,\lambda)
=\frac{2(1-\beta)(1-|\lambda|^2)(\cos\gamma)e^{-i\gamma}az^2}
{(1-\lambda z)(1+(\overline{\lambda}a-\lambda)z-az^2)},
$$
$$H_{a, \lambda}(z)+B(z,\lambda) =
\frac{2(1-\beta)(1-|\lambda|^2)(\cos\gamma)e^{-i\gamma}z}{(z-\overline{\lambda})
(1+(\overline{\lambda}a-\lambda)z-az^2)}
$$
and hence
\begin{eqnarray*}
H_{a, \lambda}(z) -c(z,\lambda) & = &
H_{a, \lambda}(z)-\frac{A(z,\lambda)+|z|^2|\tau(z,\lambda)|^2 B(z,\lambda)}
{1-|z|^2|\tau(z,\lambda)|^2}
\\
 &  =  &
\frac{1}{1-|z|^2 |\tau(z,\lambda)|^2}
\left\{\frac{}{}\left(H_{a,\lambda}(z)-A(z,\lambda)\right)\right.\\
& & \qquad \left . \frac{}{} -|z|^2|\tau(z,\lambda)|^2
\left(H_{a,\lambda}(z)+B(z,\lambda)\right)\right\}
\\
& = & \frac{2(1-\beta)(1-|\lambda|^2)(\cos\gamma)e^{-i\gamma}az^2
[\overline{1+(\overline{\lambda}a-\lambda)z-az^2}]}
{(1-|z|^2)(1+|z|^2-2{\rm Re\,}(\lambda z))
(1+(\overline{\lambda}a-\lambda)z-az^2)}\\
& = & r(z,\lambda)\frac{e^{-i\gamma}az^2}{|z|^2} \left
(\frac{|1+(\overline{\lambda}a-\lambda)z-az^2|^2}
{(1+(\overline{\lambda}a-\lambda)z-az^2)^2}\right ).
\end{eqnarray*}
Now by substituting $a=e^{i\theta}$ we easily see that
\begin{eqnarray*}
H_{e^{i\theta}, \lambda}(z) -c(z,\lambda)
& = & r(z,\lambda)\frac{e^{-i\gamma}e^{i\theta}z^2}{|z|^2} \left
(\frac{|1+(\overline{\lambda}e^{i\theta}-\lambda)z-e^{i\theta}z^2|^2}
{(1+(\overline{\lambda}e^{i\theta}-\lambda)z-e^{i\theta}
z^2)^2}\right ).
\end{eqnarray*}
For $G(z)$ as in Lemma \ref{pvdev-10-lem01}, we get that
\be\label{pvdev-10-eq15}
H_{e^{i\theta}, \lambda}(z)- c(z,\lambda) =r(z,\lambda)e^{-i\gamma}\frac{G'(z)}{|G'(z)|}
\ee
and  there exists a starlike univalent
function $G_0$ in $\mathbb{D}$ such that
$G=(1/3)e^{i\theta}G^3_0$ and $G_0(0)= G'_0(0)-1=0$.
As  $G_0$ is starlike, for any $z_0\in\mathbb{D}\setminus\{0\}$ the linear
segment joining $0$ and $G_0(z_0)$ entirely lies in
$G_0(\mathbb{D})$. Now, we define $\gamma_0$ by
\be\label{pvdev-10-eq16}
\gamma_0:\,z(t)=G_0^{-1}(tG_0(z_0)),
\quad 0\leq t \leq 1.
\ee
Since $G(z(t))=(1/3)e^{i\theta}(G_0(z(t)))^3=
(1/3)e^{i\theta}(tG_0(z_0))^3=t^3G(z_0)$, we have
\be\label{pvdev-10-eq17}
G'(z(t))z'(t)=3t^2G(z_0),\quad t\in [0,1].
\ee
Using  (\ref{pvdev-10-eq17}) and (\ref{pvdev-10-eq15}) we have
\begin{eqnarray}
\label{pvdev-10-eq18}
\int_0^{z_0}H_{e^{i\theta},\lambda}(\zeta)\,d\zeta &-& C(\lambda,\gamma_0)\\
& = & \int_0^1\left\{H_{e^{i\theta}, \lambda}(z(t))-c(z(t),\lambda)\right\}z'(t)\,dt \nonumber\\
& = & e^{-i\gamma}\int_0^1 r(z(t),\lambda)
\frac{G'(z(t))z'(t)}{|G'(z(t))z'(t)|}|z'(t)|\,dt\nonumber \\
& = & e^{-i\gamma}\frac{G(z_0)}{|G(z_0)|}\int_0^1 r(z(t),\lambda)|z'(t)|\,dt
\nonumber \\
& = & e^{-i\gamma} \frac{G(z_0)}{|G(z_0)|}R(\lambda, \gamma_0)\nonumber ,
\end{eqnarray}
where $C(\lambda,\gamma_0)$ and $R(\lambda,\gamma_0)$ are defined as in
Corollary \ref{pvdev-10-cor1}. Thus, we have
$$\int_0^{z_0}H_{e^{i\theta},\lambda}(\zeta)\,d\zeta\in\partial{\overline{\mathbb{D}}}
(C(\lambda, \gamma_0),R(\lambda, \gamma_0).
$$
Also, from Corollary \ref{pvdev-10-cor1}, we have
$$\int_0^{z_0}H_{e^{i\theta},\lambda}(\zeta)\,d\zeta\in V_{\mathcal{P}}(z_0,\lambda)
\subset\overline{\mathbb{D}}
(C(\lambda,\gamma_0),R(\lambda, \gamma_0).
$$
Hence, we  conclude that $\int_0^{z_0}H_{e^{i\theta},\lambda}(\zeta)\,d\zeta\in
\partial V_{\mathcal{P}}(z_0,\lambda)$.

Finally, we prove the uniqueness of the curve. Suppose that
$$\int_0^{z_0} P(\zeta)\,d\zeta=\int_0^{z_0}
H_{e^{i\theta},\lambda}(\zeta)\,d\zeta
$$
for some $P\in {\mathcal P}(\lambda)$ and $\theta\in (-\pi, \pi]$. We introduce
$$h(t)=e^{i\gamma}\frac{\overline{G(z_0)}}{|G(z_0)|}
\left\{P(z(t))-c(z(t),\lambda)\right\}z'(t),
$$
where $\gamma_0:\,z(t)$, $ 0\leq t \leq 1$, is given by
(\ref{pvdev-10-eq16}). Then, $h(t)$ is continuous function in
$[0,1]$ and satisfies
$$|h(t)|\leq r(z(t),\lambda)|z'(t)|.
$$
Furthermore, from (\ref{pvdev-10-eq18}) we have
\begin{eqnarray*}
\int_0^1 {\rm Re} \,h(t)\,dt
& = & \int_0^1 {\rm Re} \, \left\{e^{i\gamma}\frac{\overline{G(z_0)}}{|G(z_0)|}
\left\{P(z(t))-c(z(t),\lambda)\right\}z'(t)\right\}dt\\
& = & {\rm Re} \,\left\{e^{i\gamma}\frac{\overline{G(z_0)}}{|G(z_0)|}
\left\{\int_0^{z_0} H_{e^{i\theta},\lambda}(\zeta)\,d\zeta-C(\lambda,\gamma_0)
\right\}\right\}\nonumber\\
& = & \int_0^1 r(z(t),\lambda)|z'(t)|\,dt\nonumber.
\end{eqnarray*}
Thus, we have
$$h(t)= r(z(t),\lambda)|z'(t)| ~\mbox{ for all $t\in [0,1]$.}
$$
From (\ref{pvdev-10-eq15}) and (\ref{pvdev-10-eq17}), it follows
that
$$\int_0^{z_0}P(\zeta)\,d\zeta= \int_0^{z_0}H_{e^{i\theta},\lambda}(\zeta)\,d\zeta
~\mbox{ on $\gamma_0$.}
$$
In view of the identity theorem for analytic
functions, we see that it holds for all $z_0\in \ID$,
and hence, by the normalization, $P=H_{e^{i\theta},\lambda}$ in $\mathbb{D}$.
\epf

\section{Main Theorem}\label{sec-05}

\bthm \label{pvdev-10-th1} For $\lambda\in\mathbb{D}$ and
$z_0\in\mathbb{D}\setminus \{0\}$, the boundary
$\partial{V_{\mathcal{P}}(z_0,\lambda)}$ is the Jordan curve given by
\begin{eqnarray*}
(-\pi,\pi]\ni \theta & \mapsto & \int_0^{z_0} H_{e^{i\theta},\lambda}(\zeta)\,d\zeta \\
& = & \int_0^{z_0}\frac{1+[2(1-\beta)(\cos\gamma)e^{-i\gamma}-1]\delta
(e^{i\theta}\zeta, \lambda)\zeta}
{1-\delta(e^{i\theta}\zeta, \lambda)\zeta} \,d\zeta.
\end{eqnarray*}
If  $\int_0^{z_0} P(\zeta)\,d\zeta= \int_0^{z_0} H_{e^{i\theta},\lambda}(\zeta)\,d\zeta$
for some $P\in{\mathcal P}(\lambda)$ and $\theta\in (-\pi,\pi]$, then
$P(z)=H_{e^{i\theta},\lambda}(z)$,
where $\delta(z,\lambda)$ is defined by {\rm (\ref{pvdev-10-delta})}.
\ethm
\bpf We need to prove that the closed curve
\be\label{pvdev-10-curve}
(-\pi,\pi]\ni \theta \mapsto \int_0^{z_0}H_{e^{i\theta},\lambda}(\zeta)\,d\zeta
\ee
is simple. Suppose that
$$\int_0^{z_0}H_{e^{i\theta_1},\lambda}(\zeta)\,d\zeta=\int_0^{z_0}
H_{e^{i\theta_2},\lambda}(\zeta)\,d\zeta
$$
for some $\theta_1,\theta_2\in(-\pi,\pi]$ with $\theta_1\neq\theta_2$.
Then, from Proposition \ref{pvdev-10-pro2}, we have
\be\label{pvdev-10-eq21}
H_{e^{i\theta_1},\lambda}= H_{e^{i\theta_2},\lambda}.
\ee
From (\ref{pvdev-10-eq5}) and (\ref{pvdev-10-eq11}) we obtain the
following identity
\be\label{pvdev-10-eq22}
\tau\left(\frac{\omega_{H_{e^{i\theta},\lambda}}}{z},\lambda\right)=
\frac{e^{i\theta}z(1-\overline{\lambda}^2)+\lambda-\overline{\lambda}}
{e^{i\theta}z(\overline{\lambda}-\lambda)+1-{\lambda}^2}.
\ee
From (\ref{pvdev-10-eq21}) and (\ref{pvdev-10-eq22}) we
have
\be\label{pvdev-10-eq23}
\frac{e^{i\theta_1}z(1-\overline{\lambda}^2)+\lambda-\overline{\lambda}}
{e^{i\theta_1}z(\overline{\lambda}-\lambda)+1-{\lambda}^2}=\frac{e^{i\theta_2}
z(1-\overline{\lambda}^2)+\lambda-\overline{\lambda}}
{e^{i\theta_2}z(\overline{\lambda}-\lambda)+1-{\lambda}^2}.
\ee
A simplification of (\ref{pvdev-10-eq23}) gives
$$e^{i\theta_1}z=e^{i\theta_2}z
$$
which is a contradiction to the choice of $\theta_1$ and $\theta_2$.
Thus, the curve must be simple.

Since $V_{\mathcal{P}}(z_0,\lambda)$ is a compact convex subset of $\mathbb{C}$
and has nonempty interior, the boundary $\partial V_{\mathcal{P}}(z_0,\lambda)$
is a simple closed curve. From Proposition \ref{pvdev-10-pro1}, the
curve $\partial V_{\mathcal{P}}(z_0,\lambda)$ contains the curve
(\ref{pvdev-10-curve}).
Recall the fact that a simple closed curve cannot contain any simple closed curve
other than itself. Thus, $\partial V_{\mathcal{P}}(z_0,\lambda)$ is given by
(\ref{pvdev-10-curve}).
\epf

\br\label{rem1}
The integral in (\ref{pvdev-10-curve}) can be simplified as follows:
Set $b = {\rm Im}(\overline{\lambda}e^{i\theta/2}) \in {\mathbb
R}$. Then a computation shows that
\beqq
1+(\overline{\lambda}e^{i\theta}-\lambda)z-e^{i\theta}z^2
&=& (1- z/z_1)(1-z/z_2 ),
\eeqq
where
$$z_1 =e^{-i\theta /2}(ib+\sqrt{1-b^2}) ~\mbox{ and } ~
z_2=e^{-i\theta /2}(ib-\sqrt{1-b^2}).
$$
From (\ref{pvdev-10-extremal})  and (\ref{pvdev-10-delta}) we have
\be\label{pvdev-8-eq001}
 H_{e^{i\theta},\lambda}(z)  =
1+2(1-\beta)e^{-i\gamma}\cos\gamma \left(\frac{(e^{i\theta}z+\lambda)z}
{1+(\overline{\lambda}e^{i\theta}-\lambda)z-e^{i\theta}z^2}\right).
\ee
Since
$$
\frac{(e^{i\theta}z+\lambda)z}{1+(\overline{\lambda}e^{i\theta}-\lambda)z-e^{i\theta}z^2}
=-1-\frac{e^{-i\theta}}{z_1-z_2}
\left(\frac{1+\overline{\lambda}e^{i\theta}z_1}{z-z_1}
-\frac{1+\overline{\lambda}e^{i\theta}z_2}{z-z_2} \right),
$$
the equation (\ref{pvdev-8-eq001}) becomes
\begin{eqnarray*}
H_{e^{i\theta},\lambda}(z) & = & 1-2(1-\beta)e^{-i\gamma}\cos\gamma\\
& & -\frac{2e^{-i\gamma}(1-\beta)e^{-i\theta}\cos\gamma}{z_1-z_2}
\left(\frac{1+\overline{\lambda}e^{i\theta}z_1}{z-z_1}
-\frac{1+\overline{\lambda}e^{i\theta}z_2}{z-z_2} \right).
\end{eqnarray*}
By integrating on both sides from $0$ to $z_0$, we can easily obtain
the following representation:
\begin{eqnarray*}
& & \int_0^{z_0} H_{e^{i\theta},\lambda}(\zeta)\,d\zeta=(1-2(1-\beta)e^{-i\gamma}\cos\gamma)z_0 +\\
& & K(\gamma,\beta,\theta,b)
\left[\left(1+\overline{\lambda}e^{i\theta/2}(-\sqrt{1-b^2}+ib)\right)
\log\left(1+\frac{e^{i\theta/2}z_0}{\sqrt{1-b^2}-ib}\right)\right.\\
&& \quad\left.
-\left(1+\overline{\lambda}e^{i\theta/2}(\sqrt{1-b^2}+ib)\right)
\log\left(1-\frac{e^{i\theta/2}z_0}{\sqrt{1-b^2}+ib}\right)
\right],
\end{eqnarray*}
where
$$K(\gamma,\beta,\theta,b)=\frac{e^{-i\gamma}(1-\beta)e^{-i\theta/2}\cos\gamma }{\sqrt{1-b^2}}.
$$

\er

\vspace{8pt}

For $\lambda=0$,  Theorem \ref{pvdev-10-th1} takes the following simple form.

\bcor
For $z_0\in\mathbb{D}\setminus \{0\}$ and $\lambda=0$ the boundary
$\partial{V_{\mathcal{P}}(z_0,0)}$ is the Jordan curve given by
\begin{eqnarray*}
& & (-\pi,\pi]\ni \theta  \mapsto
\int_0^{z_0} H_{e^{i\theta},0}(\zeta)\,d\zeta \\
& = & (1-2(1-\beta)e^{-i\gamma}\cos\gamma)z_0
+e^{-i\gamma}(1-\beta)e^{-i\theta/2}\cos\gamma \log\left(\frac{1+e^{i\theta/2} z_0}{1-e^{i\theta/2} z_0} \right)
\end{eqnarray*}
If  $\int_0^{z_0} P(\zeta)\,d\zeta= \int_0^{z_0} H_{e^{i\theta},0}(\zeta)\,d\zeta$
for some $P\in{\mathcal P}(0)$ and $\theta\in (-\pi,\pi]$, then
$P(z)=H_{e^{i\theta},0}(z)$.
\ecor

\section{Some special cases} \label{sec-06}
\subsection{The class $\mathcal{R}_{\beta}$}
In order to discuss a special situation, we consider $P=f'$ and $\gamma=0$
in the class $\mathcal{P}_{\gamma, \beta}$. Thus,  $\mathcal{P}_{\gamma, \beta}$
reduces to $\mathcal{R}_{\beta}$, where
$$\mathcal{R}_{\beta}=\{f\in\mathcal{A}\colon {\rm Re\,}f'(z)>\beta
\quad\mbox{ in } \mathbb{D}\}.
$$
Then $\mathcal{R}_\beta \subset \mathcal{S}$ for $0\leq\beta<1$.
As with $\mathcal{P}(\lambda )$,  for $\lambda\in\overline{\mathbb{D}}$ and
$z_0\in\mathbb{D}$ being fixed,  we define
\beqq
\mathcal{R}(\lambda)  & = &
\left\{\frac{}{} f\in\mathcal{R}_{\beta} :\,
f''(0)=2(1-\beta)\lambda
\right\}\\
V_{\mathcal{R}}(z_0,\lambda) & = & \left\{\frac{}{} f(z_0)\colon
f\in\mathcal{R}(\lambda)\right\}.
\eeqq
We remark that if $f \in \mathcal{R}(\lambda)$, then it is necessary that $f'''(0)$ satisfies
the condition
$$f'''(0)=4(1-\beta)[(1-|\lambda|^2)a+ {\lambda}^2]
$$
for some $a\in\overline{\mathbb{D}}$.

For $P=f'$, a computation shows that the extremal function $H_{e^{i\theta},\lambda}(z)$
for the class $\mathcal{R}(\lambda)$ takes the form
\begin{eqnarray*}
H_{e^{i\theta},\lambda}(z)=z_0+2(1-\beta)\int_0^{z_0}\frac{(e^{i\theta}\zeta+\lambda)\zeta}
{1+\overline{\lambda}e^{i\theta}\zeta-(e^{i\theta}\zeta+\lambda)\zeta} \,d\zeta.
\end{eqnarray*}

It is not difficult to obtain the following result which is the analog of
Theorem \ref{pvdev-10-th1} for the class $\mathcal{R}(\lambda)$.

\bcor
For $\lambda\in\mathbb{D}$ and $z_0\in\mathbb{D}\setminus\{0\}$, the boundary
$\partial V_{\mathcal{R}}(z_0,\lambda)$ is the Jordan curve given by
\begin{eqnarray*}
(-\pi,\pi]\ni \theta \mapsto H_{e^{i\theta},\lambda}(z_0)= z_0+2(1-\beta)\int_0^{z_0}
\frac{(e^{i\theta}\zeta+\lambda)\zeta}
{1+\overline{\lambda}e^{i\theta}\zeta-(e^{i\theta}\zeta+\lambda)\zeta} \,d\zeta.
\end{eqnarray*}
If  $f(z_0)=H_{e^{i\theta},\lambda}(z_0)$ for some
$f\in\mathcal{R}(\lambda)$ and $\theta\in (-\pi,\pi]$, then
$f(z)=H_{e^{i\theta},\lambda}(z)$.
\ecor

For $0\leq\beta<1$ and $\lambda =0$, set
$$\mathcal{R}(0 )=\{f\in\mathcal{A}\colon f''(0)=0 ~\mbox{ and }
{\rm Re\,}f'(z)>\beta ~\mbox{ in } \mathbb{D}\} \subset \mathcal{R}_{\beta} .
$$
In  particular, the choices $\gamma=0$ and $P(z)=f'(z)$ in Corollary
\ref{pvdev-10-cor01} give the following: if $f\in\mathcal{R}(0) \subset \mathcal{R}_{\beta} $
for some $0\leq\beta<1/2$, then by (\ref{pvdev-10-eq24}), one has
$$|f'(z)|\leq \frac{1+(1-2\beta)|z|^4 +2(1-\beta)|z|^2}{1-|z|^4}
= \frac{1+(1-2\beta)|z|^2}{1-|z|^2}, \quad z\in \ID ,
$$
so that
$$\sup_{z\in\mathbb{D}}(1-|z|^2)|f'(z)|\leq 2(1-\beta).
$$
Equality holds for
$$f(z)=\beta z+\frac{(1-\beta)}{2}\log\left(\frac{1+z}{1-z}\right), \quad z\in\mathbb{D}.
$$
 \subsection{The class ${\mathcal F}(\alpha, \beta)$}

For a complex number $\alpha\in\mathbb{C}$ satisfying ${\rm Re\,}{\alpha}>0$
and $\beta\in\mathbb{R}$ with $\beta<1$, let ${\mathcal F}(\alpha, \beta)$ denote
the class of functions $f\in {\mathcal A}$ satisfying
\be\label{pvdev-8-eq1}
f'(z)+\alpha zf''(z)\prec \frac{1+(1-2\beta)z}{1-z},\quad z\in\mathbb{D},
\ee
where $\prec$ denote the usual subordination \cite{Miller and Mocanu-book}.
In \cite{P5} conditions on $\alpha$ and $\beta$ for which
$${\mathcal F}(\alpha, \beta)\subset {\mathcal S}^*
$$
have been established (including for certain  complex values of $\alpha$)
and in \cite{FR1} it has been shown that
${\mathcal F}(\alpha, \beta)\subset {\mathcal S}^*$
if $\alpha\geq 1/3$ and  $\beta \geq \beta_0(\alpha)$, where
$$\beta_0(\alpha)=\frac{-\frac{1}{\alpha}
\int_0^1t^{\frac{1}{\alpha}-1}\left(\frac{1+t}{1-t}\right)\,dt}
{1-\frac{1}{\alpha}\int_0^1t^{\frac{1}{\alpha}-1}\left(\frac{1+t}{1-t}\right)\,dt}.
$$
This is indeed a reformulated version of a theorem from \cite{FR1} and the inclusion
is sharp in the following  sense: for $\beta<\beta_0(\alpha)$ the functions
in $\mathcal{F}(\alpha,\beta)$ are not even univalent in $\mathbb{D}$.
For an extension of this  inclusion result,  we refer to \cite{PR1,PR2}.

Now, we present an alternative  representation for functions in
$\mathcal{F}(\alpha, \beta)$. If $f\in\mathcal{F}(\alpha, \beta)$, then
(\ref{pvdev-8-eq1}) is equivalent to
$$
\frac{f(z)}{z}* \left(1+\sum_{n=2}^{\infty} n(1+(n-1)\alpha)z^{n-1}\right)\prec
1+2(1-\beta)\frac{z}{1-z},
$$
where $*$ denotes the Hadamard product (or convolution) of two analytic
functions in $\ID$ represented by power series about the origin.
By a well-known convolution theorem (cf. \cite{rush-stan85}) this gives
$$\frac{f(z)}{z}\prec
\beta+(1-\beta)\left[1+\frac{2}{\alpha}
\sum_{n=1}^{\infty}\frac{z^n}{(n+1)(n+1/{\alpha})}\right]
$$
and a computation shows that
$$
\frac{f(z)}{z}\prec\left\{\ba{ll}
\ds
\beta+(1-\beta)\left[1-\frac{2}{1-\alpha}\left(\frac{\log(1-z)}{z}
+1+\int_0^1 t^{1/\alpha}\frac{z}{1-tz}\,dt \right) \right]\\[6mm]
& \hspace*{-2cm} \mbox{if $\alpha\neq 1$} \\ [6mm]
\ds
\beta+(1-\beta)\left[1+2z\int_0^1 \frac{t\log(1/t)}{1-tz}\,dt  \right]
&  \hspace*{-2cm} \mbox{if $\alpha=1$.}
\ea \right .
$$

The definition of subordination gives the following representation of functions in
$\mathcal{F}(\alpha, \beta)$:
$$
f(z)=\left\{\ba{ll}
\ds
z-\frac{2(1-\beta)z}{1-\alpha}\left\{1+\frac{1}{\omega(z)}\log(1-\omega(z))+\omega(z)
\int_0^1\frac{t^{1/\alpha}}{1-t\omega(z)}\,dt\right\}\\[6mm]
&  \hspace*{-2cm}\mbox{if $\alpha\neq 1$} \\ [6mm]
\ds z+2(1-\beta)z\omega(z)\int_0^1\frac{t\log (1/t)}{1-t\omega(z)}\,dt
&  \hspace*{-2cm} \mbox{if $\alpha=1$}
\ea \right .
$$
for $z\in\mathbb{D}$, and for some  $\omega \in {\mathcal B}_0$.

If $f\in{\mathcal F}(\alpha, \beta)$, then according to the Herglotz representation
there exists a unique positive unit measure $\mu$ on $(-\pi,\pi]$ such that
\begin{equation*}\label{pvdev-8-herglotz}
f'(z)+\alpha zf''(z)=\int_{-\pi}^{\pi}\frac{1+(1-2\beta)ze^{-it}}{1-ze^{-it}}\,d\mu(t),
\end{equation*}
or equivalently
$$
\frac{f(z)}{z}= \left[1+\frac{1}{\alpha}\sum_{n=1}^{\infty}\frac{z^n}{(n+1)(n+1/\alpha)}
\right]\ast\int_{-\pi}^{\pi}\frac{1+(1-2\beta)ze^{-it}}{1-ze^{-it}}\,d\mu(t).
$$
A simplification of the last equality gives the following
representation of functions in the class $\mathcal{F}(\alpha,\beta)$:
$$
f(z)=\left\{\ba{ll}
\ds
\frac{z}{1-\alpha}\int_0^1\int_{-\pi}^{\pi} \left(1-s^{\frac{1}{\alpha}-1}\right)
\left(\frac{1+(1-2\beta)sze^{-it}}{1-sze^{-it}} \right) d\mu(t)\,ds \\[6mm]
& \hspace{-2cm}\mbox{if $\alpha\neq 1$} \\ [6mm]
\ds z+2(1-\beta)z\int_0^1\int_{-\pi}^{\pi}\left (\log ( 1/s)\right)
\left(\frac{sze^{-it}}{1-sze^{-it}} \right) d\mu(t)\,ds\\[6mm]
& \hspace{-2cm} \mbox{if $\alpha=1$.}
\ea \right .
$$

To state our special case in precise form, we for convenience let $\gamma=0$, and let $P$ be defined by
$$P(z)=f'(z)+\alpha zf''(z), \quad f\in\mathcal{F}(\alpha,\beta)
$$
so that
$$P'(0)=(1+\alpha)f''(0) ~\mbox{ and }~ P''(0)=(1+2\alpha)f''(0).
$$
In view of these observations, the analog of the sets $\mathcal{P}(\lambda)$ and
$V_{\mathcal{P}}(z_0,\lambda)$ will be as follows:
\beqq
\mathcal{G}(\lambda)&=&\left\{f\in{\mathcal F}(\alpha, \beta) :\,f''(0)
=2\left(\frac{1-\beta}{1+\alpha}\right)\lambda \right\}
\eeqq
and
\begin{eqnarray*}
V_{\mathcal{G}}(z_0,\lambda) &  = &  \{(1-\alpha)f(z_0)+\alpha z_0f'(z_0):\,
f\in {\mathcal G}(\lambda)\},
\end{eqnarray*}
where $0\leq \beta <1$. We observe that for functions in $\mathcal{G}(\lambda)$, $f'''(0)$ will be
of the form
$$f'''(0)= 4\big((1-|\lambda|^2)a+{\lambda}^2\big)\left(\frac{1-\beta}{1+2\alpha}\right)
$$
for some $a\in \overline{\ID}$.

With $P(z)=f'(z)+\alpha zf''(z)$, the corresponding extremal function
$F_{e^{i\theta},\lambda}(z)$ for ${\mathcal G}(\lambda)$ can be computed and this is given by
$$(1-\alpha)F_{e^{i\theta},\lambda}(z)+\alpha
zF'_{e^{i\theta},\lambda}(z) = \int_0^{z}
\frac{1+(1-2\beta)\delta(a\zeta, \lambda)\zeta}{1-\delta(a\zeta, \lambda)\zeta}\,d\zeta,
$$
where $\delta(z,\lambda)$ is defined by (\ref{pvdev-10-delta}).
In this setting, Proposition \ref{pvdev-10-pro1} (for $\gamma=0$) takes the following form:

\bprop\label{pvdev-8-pro1}
For $f\in {\mathcal G}(\lambda)$ and $\lambda\in\mathbb{D}$, we have
$$
\left|f'(z)+\alpha zf''(z)-c(z,\lambda)\right|\leq r(z,\lambda),
\quad z\in\mathbb{D},
$$
where
\begin{eqnarray*}
c(z,\lambda)& = &
\frac{(1+(1-2\beta)\lambda z)(1-\overline{\lambda}\overline{z})+
|z|^2(\overline{z}-\lambda)(\overline{\lambda}+(1-2\beta)z)}{(1-|z|^2)
(1+|z|^2-2{\rm Re\,}(\lambda z))},~~ \mbox{ and }\\
r(z,\lambda) & = & \frac{2(1-\beta)(1-|\lambda|^2)|z|^2}
{(1-|z|^2)(1+|z|^2-2{\rm Re\,}(\lambda z))}.
\end{eqnarray*}
For each $z\in\mathbb{D}\setminus\{0\}$, equality holds if and
only if $f=F_{e^{i\theta},\lambda}$ for some
$\theta\in\mathbb{R}$.
\eprop

Using Theorem \ref{pvdev-10-th1}, we get the following result.

\bcor \label{pvdev-8-th1} For $\lambda\in\mathbb{D}$, $z_0\in\mathbb{D}\setminus \{0\}$
and $\alpha\in\mathbb{C}$ with ${\rm Re\,}{\alpha}>0$, the boundary
$\partial{V_{\mathcal{G}} (z_0,\lambda)}$ is the Jordan curve given by
\begin{eqnarray*}
& & (-\pi,\pi]\ni \theta \mapsto
(1-\alpha)F_{e^{i\theta},\lambda}(z_0)+\alpha z_0F'_{e^{i\theta},\lambda}(z_0)= (2\beta-1)z_0\\
& &  +\frac{(1-\beta)e^{-i\theta/2}}{\sqrt{1-b^2}}
\left[\left(1+\overline{\lambda}e^{i\theta/2}(-\sqrt{1-b^2}+ib)\right)
\log\left(1+\frac{e^{i\theta/2}z_0}{\sqrt{1-b^2}-ib}\right)\right.\\
&& \quad\left.
-\left(1+\overline{\lambda}e^{i\theta/2}(\sqrt{1-b^2}+ib)\right)
\log\left(1-\frac{e^{i\theta/2}z_0}{\sqrt{1-b^2}+ib}\right)
\right],
\end{eqnarray*}
where $b = {\rm Im}(\overline{\lambda}e^{i\theta/2})$.
If  $(1-\alpha)f(z_0)+\alpha z_0f'(z_0)= (1-\alpha)F_{a,\lambda}(z_0)
+\alpha z_0F'_{a,\lambda}(z_0)$ for some
$f\in{\mathcal G}(\lambda)$ and $\theta\in (-\pi,\pi]$, then
$f(z)=F_{e^{i\theta},\lambda}(z)$.
\ecor
The proof of this corollary follows by taking  $\gamma=0$ in  Remark \ref{rem1}
and so we omit the details.

In the case of $\lambda =0$ in Corollary \ref{pvdev-8-th1},
the corresponding extremal function $F_{a,0}(z)$ can be obtained easily by solving
$$
F'_{a, 0}(z)+\alpha zF''_{a, 0}(z) = \frac{1+(1-2\beta)\delta(az,0)z}{1-\delta(az,0)z}.
$$
This may be rewritten as
$$\frac{F_{a,0}(z)}{z}*\left[1+\sum_{n=0}^{\infty}(n+1)(1+n\alpha)z^n\right]
 = 2\beta-1 +2(1-\beta)\frac{1}{1-az^2}
$$
or equivalently as
$$
\frac{F_{a,0}(z)}{z} =
\left[1+\sum_{n=1}^{\infty} 2(1-\beta) a^nz^{2n}\right]*
\left[1+\sum_{n=0}^{\infty}\frac{z^n}{(n+1)(1+n\alpha)}\right].
$$
A simple calculation gives that
$$
F_{a,0}(z)= \left\{\ba{ll}
\ds
z+\frac{(1-\beta)az^3}{(1-\alpha)}\int_0^1 \frac{t^{1/2} - t^{1/2\alpha}}{1-taz^2}\,dt
& \mbox{if $\alpha\neq 1$} \\ [6mm]
\ds
z+\frac{(1-\beta)az^3}{2}\int_0^1 \frac{t^{1/2}\log(1/t)}{1-taz^2}\,dt
& \mbox{if $\alpha=1$.}
\ea \right .
$$

\section{Geometric view of Theorem \ref{pvdev-10-th1}}\label{sec-07}

  Using Mathematica (see \cite{Ruskeepaa}),  we describe the boundary of the sets
$V_{\mathcal{P}}(z_0, \lambda)$ and $V_{\mathcal{G}}(z_0, \lambda)$.
In the  program below, ``z0 stands for $z_0$'',  ``lam for $\lambda$''  ``g for $\gamma$'' and
``b for $\beta$''.

\vspace{0.1cm}
{\tt
\begin{verbatim}
(* Geometric view the main Theorem 5.1  and Corollary 6.3 *)

Remove["Global`*"];

z0 = Random[]Exp[I*Random[Real, {-Pi, Pi}]]
lam = Random[]Exp[I*Random[Real, {-Pi, Pi}]]
g = Random[Real, {-Pi/2, Pi/2}]
b = Random[Real, {0, 1}]

Print["z0=", z0]
Print["lam=", lam]
Print["g=", g]
Print["b=", b]


Q1[b_, g_, lam_,the_] := ((1 +Conjugate[lam]Exp[I*the]*z) +
(lam*z +Exp[I*the]*z^2)(Exp[-I*g] - 2b*Cos[g])Exp[-I*g])/
((1 + ( Conjugate[lam]*Exp[I*the] - lam)*z )-Exp[I*the]*z*z);


Q2[b_, lam_,the_] := ((1 +Conjugate[lam]Exp[I*the]*z) +
(1 - 2b)(lam*z +Exp[I*the]*z^2))/
((1 + ( Conjugate[lam]*Exp[I*the] - lam)*z )-Exp[I*the]*z*z);

myf1[b_, g_, lam_, the_, z0_] :=
NIntegrate[Q1[b, g, lam, the], {z, 0, z0}];

myf2[b_, lam_, the_, z0_] := NIntegrate[Q2[b, lam, the],
                             {z, 0, z0}];

image1 = ParametricPlot[{Re[myf1[b, g,  lam, the, z0]],
Im[myf1[b, g, lam, the, z0]]}, {the, -Pi, Pi},AspectRatio ->
Automatic,TextStyle -> {FontFamily -> "Times",FontSize ->14},
AxesStyle -> {Thickness[0.0035]} ];

image2 = ParametricPlot[{Re[myf2[b,  lam, the, z0]],
Im[myf2[b, lam, the, z0]]},{the, -Pi, Pi}, AspectRatio ->
Automatic,TextStyle -> {FontFamily -> "Times",FontSize ->14},
AxesStyle -> {Thickness[0.0035]} ];

Clear[b, g, lam, z0, myf1, myf2];
\end{verbatim}
}

\vspace{0.5cm}

Following figures show the boundary of $V_{\mathcal{P}}(z_0, \lambda)$  and
 $V_{\mathcal{G}}(z_0, \lambda)$ for certain values
of  $z_0\in\mathbb{D}\setminus\{0\}$, $\lambda\in\mathbb{D}$, $0\leq\beta<1$ and
$|\gamma|<\pi/2$. Table 1 gives the list of these parameter values corresponding to Figs. 1-5.
We recall that according to  Proposition \ref{pvdev-10-pro01}
the region bounded by the curve $\partial V_{\mathcal{P}}(z_0, \lambda)$
is compact and convex.

\begin{center}
Table 1
\end{center}

\begin{center}
\begin{tabular}{|c|c|c|c|c|}
\hline
Fig. & $z_0$ & $\lambda$ & $\beta$  & $\gamma$  \\ \hline
1 & 0.335192-0.787333i & 0.0737292+0.466706i & 0.591244 & 0.383292 \\\hline
2 & -0.261209+0.926935i & -0.28588+0.307498i  & 0.700318  & -0.87825  \\\hline
3 & -0.41227-0.521734i   & -0.0875648+0.0714166i   & 0.602203   & 0.910581  \\ \hline
4 & 0.771264+0.151204i  & -0.391149-0.294747i   & 0.928608   & 1.55854  \\ \hline
5 & 0.335626+0.929093i   & 0.00010443+0.0255256i   & 0.76622   & 1.5449  \\ \hline
\end{tabular}
\end{center}

\begin{figure}[htp]
\begin{center}
\includegraphics[width=6.2cm]{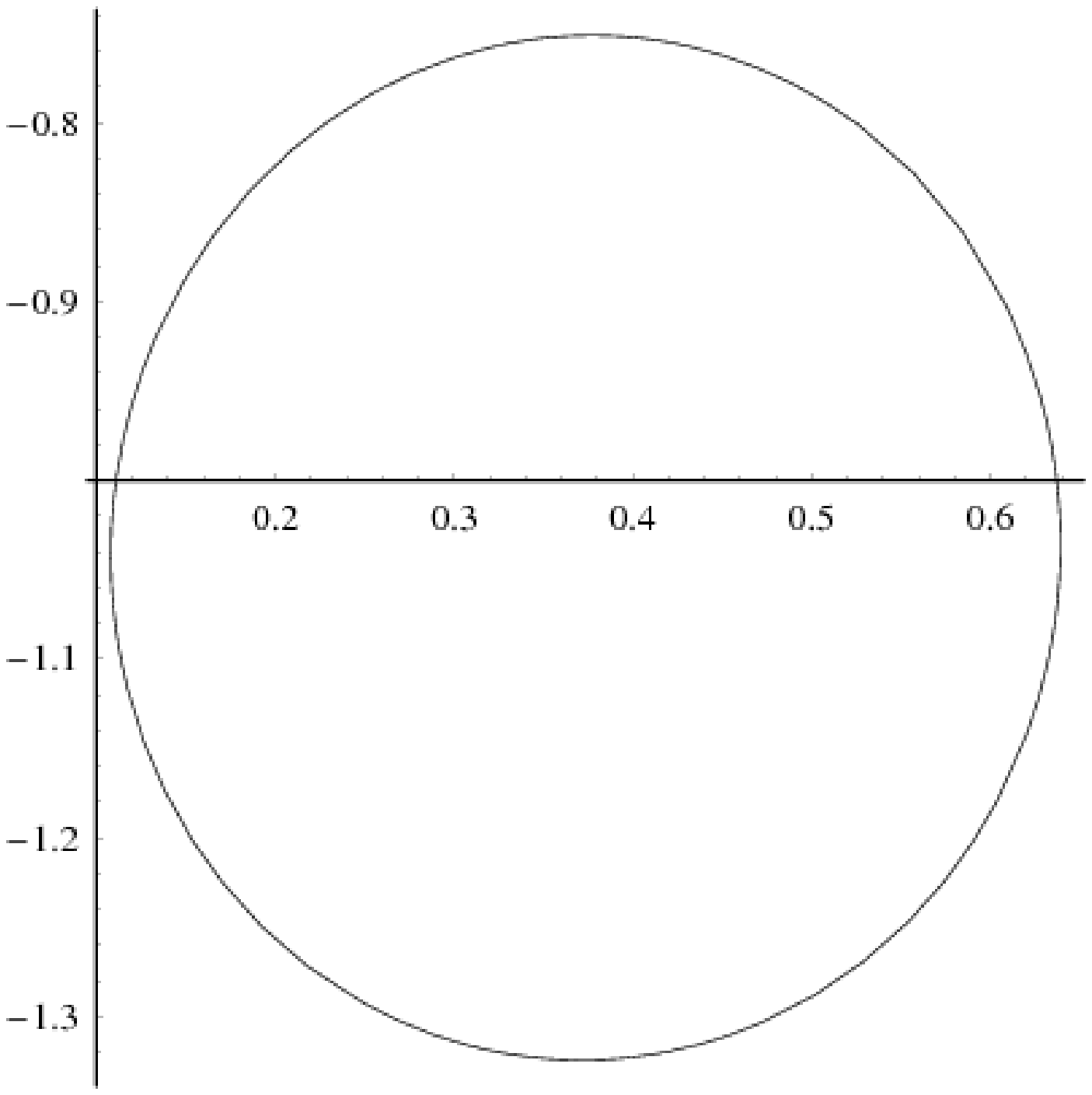}
\hspace{0.5cm}
\includegraphics[width=6.2cm]{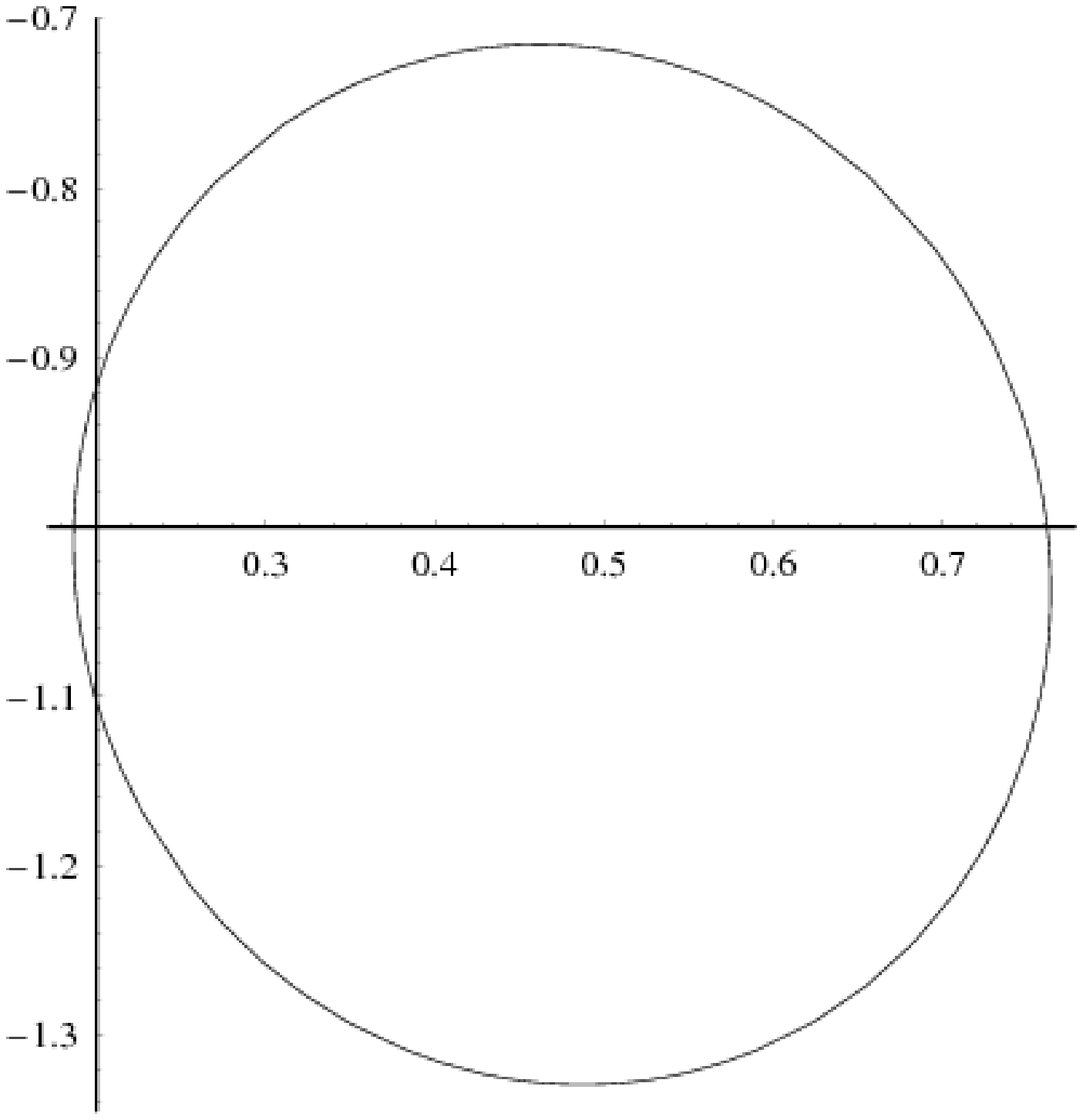}
\end{center}
\hspace{1.5cm}$\partial V_{\mathcal P}(z_0,\lambda)$ \hspace{4.5cm}  $\partial
V_{\mathcal G}(z_0,\lambda)$
\caption{$z_0 = 0.335192-0.787333i$ and $\beta =0.591244$}
\end{figure}


\begin{figure}[htp]
\begin{center}
\includegraphics[width=6.6cm]{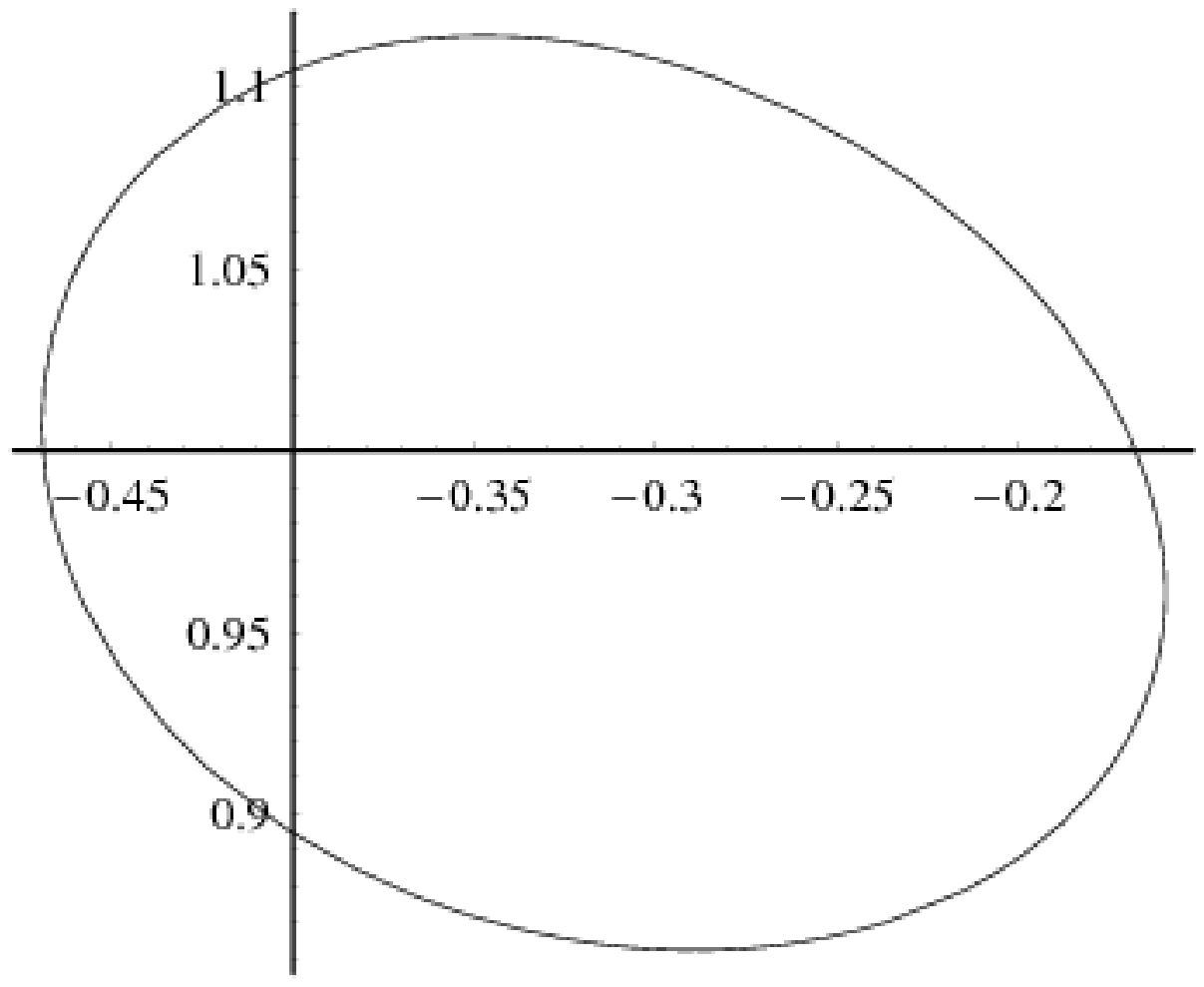}
\hspace{-0.5cm}
\includegraphics[width=6.6cm]{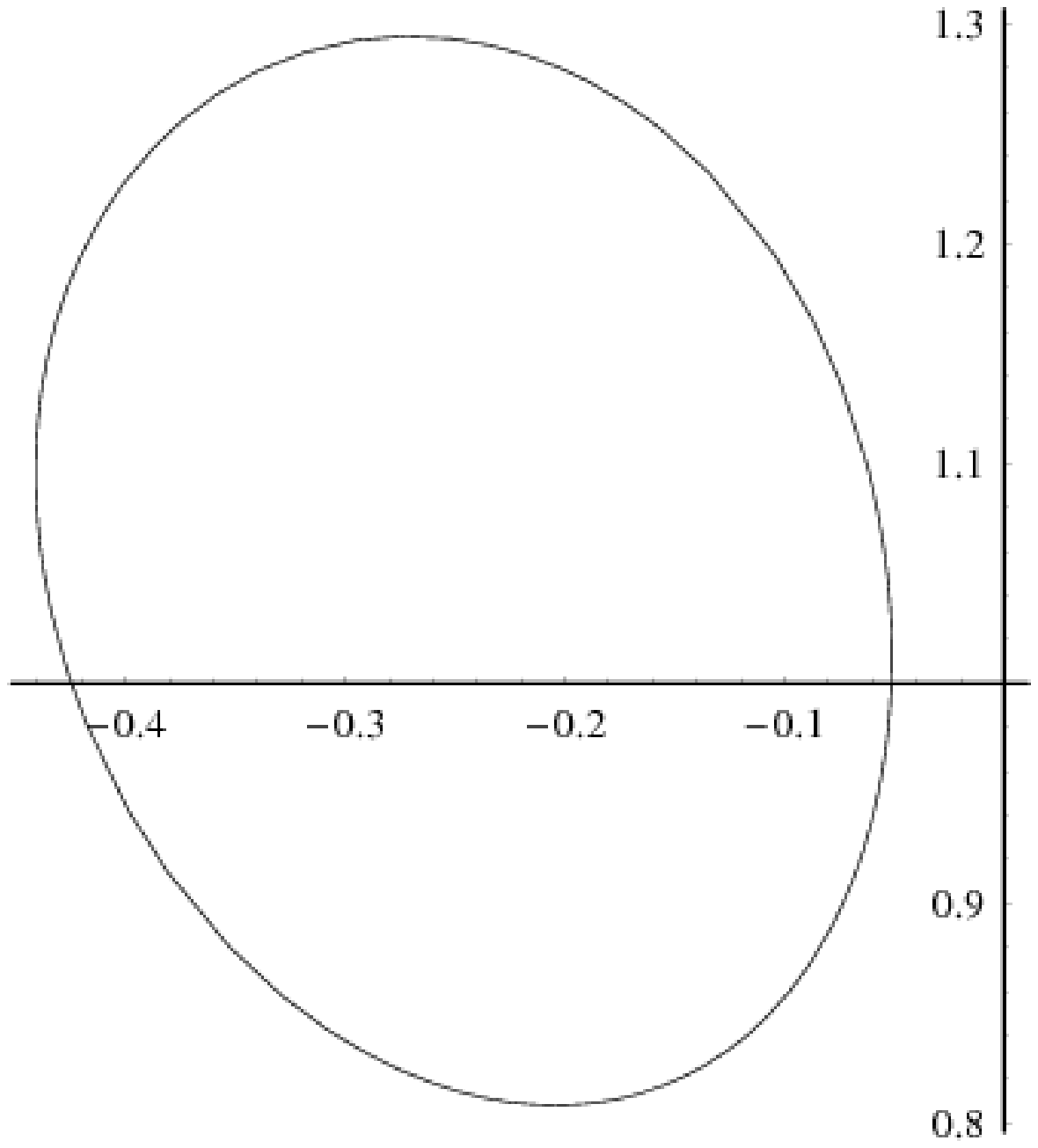}
\end{center}
\hspace{1.5cm}$\partial V_{\mathcal P}(z_0,\lambda)$ \hspace{4.5cm}  $\partial
V_{\mathcal G}(z_0,\lambda)$
\caption{$z_0 =-0.261209+0.926935i$ and $\beta= 0.700$}
\end{figure}


\begin{figure}[htp]
\begin{center}
\includegraphics[width=6.2cm]{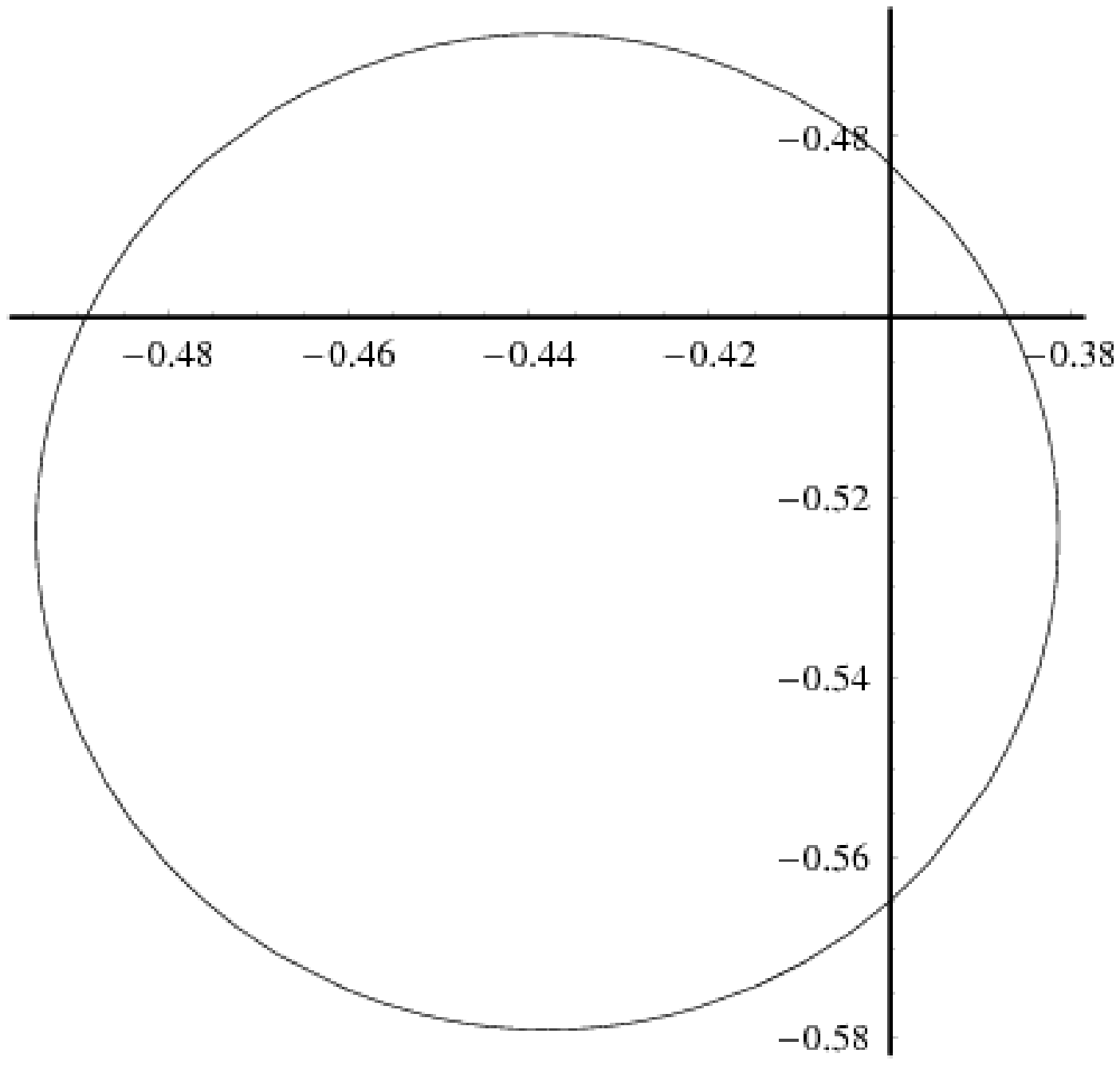}
\hspace{0.5cm}
\includegraphics[width=6.2cm]{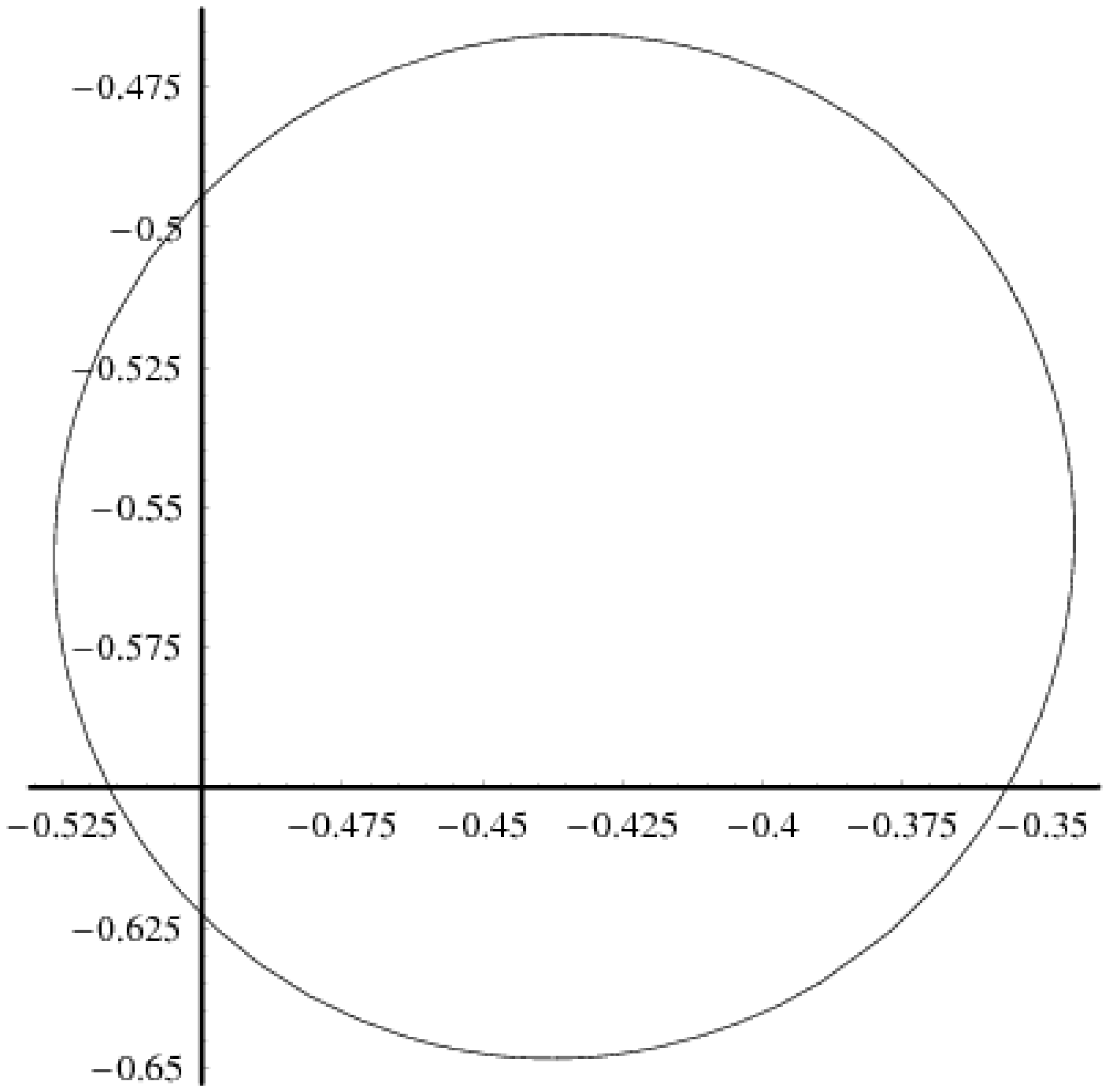}
\end{center}
\hspace{1.5cm}$\partial V_{\mathcal P}(z_0,\lambda)$ \hspace{4.5cm}  $\partial
V_{\mathcal G}(z_0,\lambda)$
\caption{$z_0 =-0.41227-0.521734i$ and $\beta =0.602203$}
\end{figure}


\begin{figure}[htp]
\begin{center}
\includegraphics[width=6.6cm]{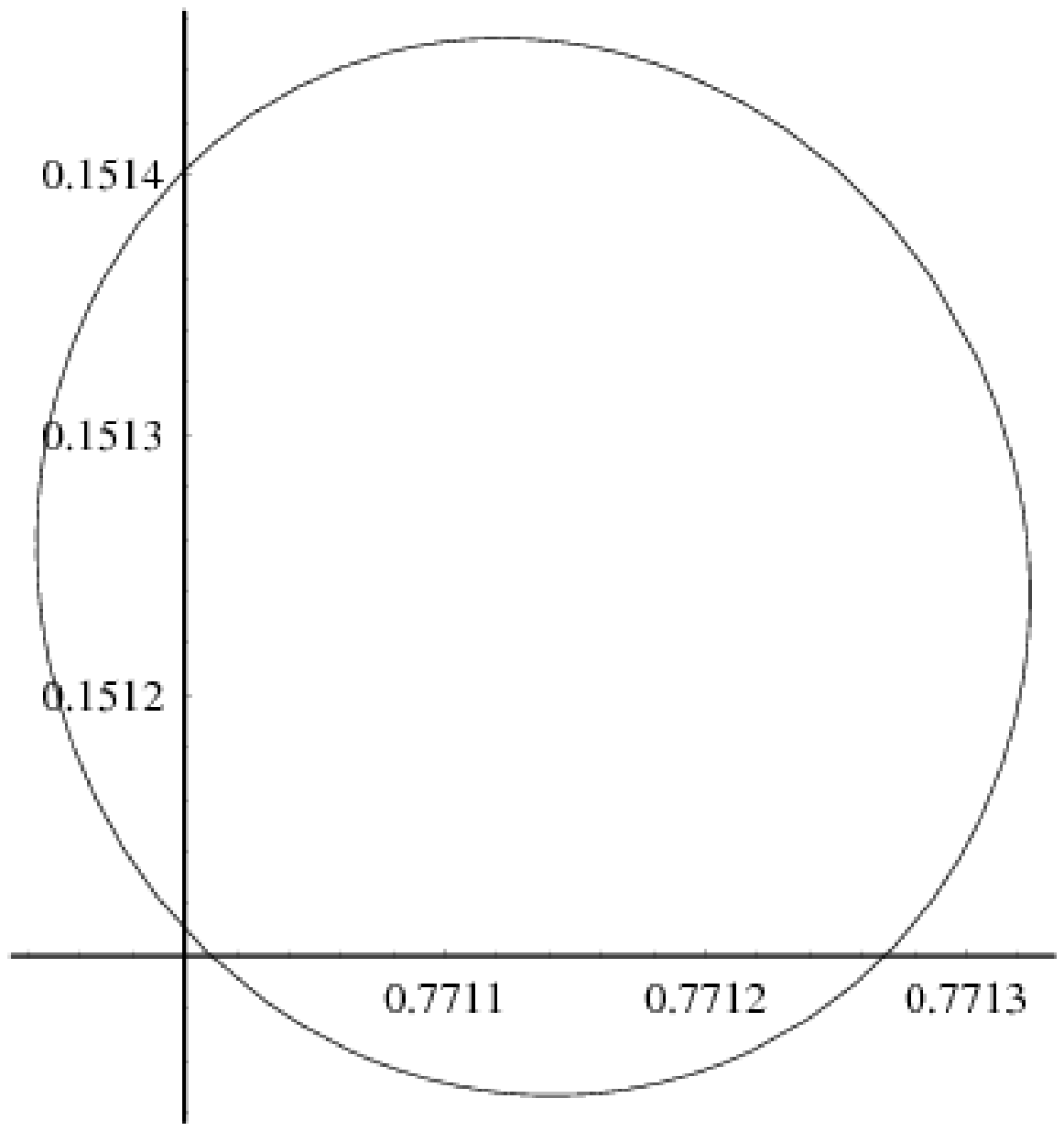}
\hspace{-0.5cm}
\includegraphics[width=6.6cm]{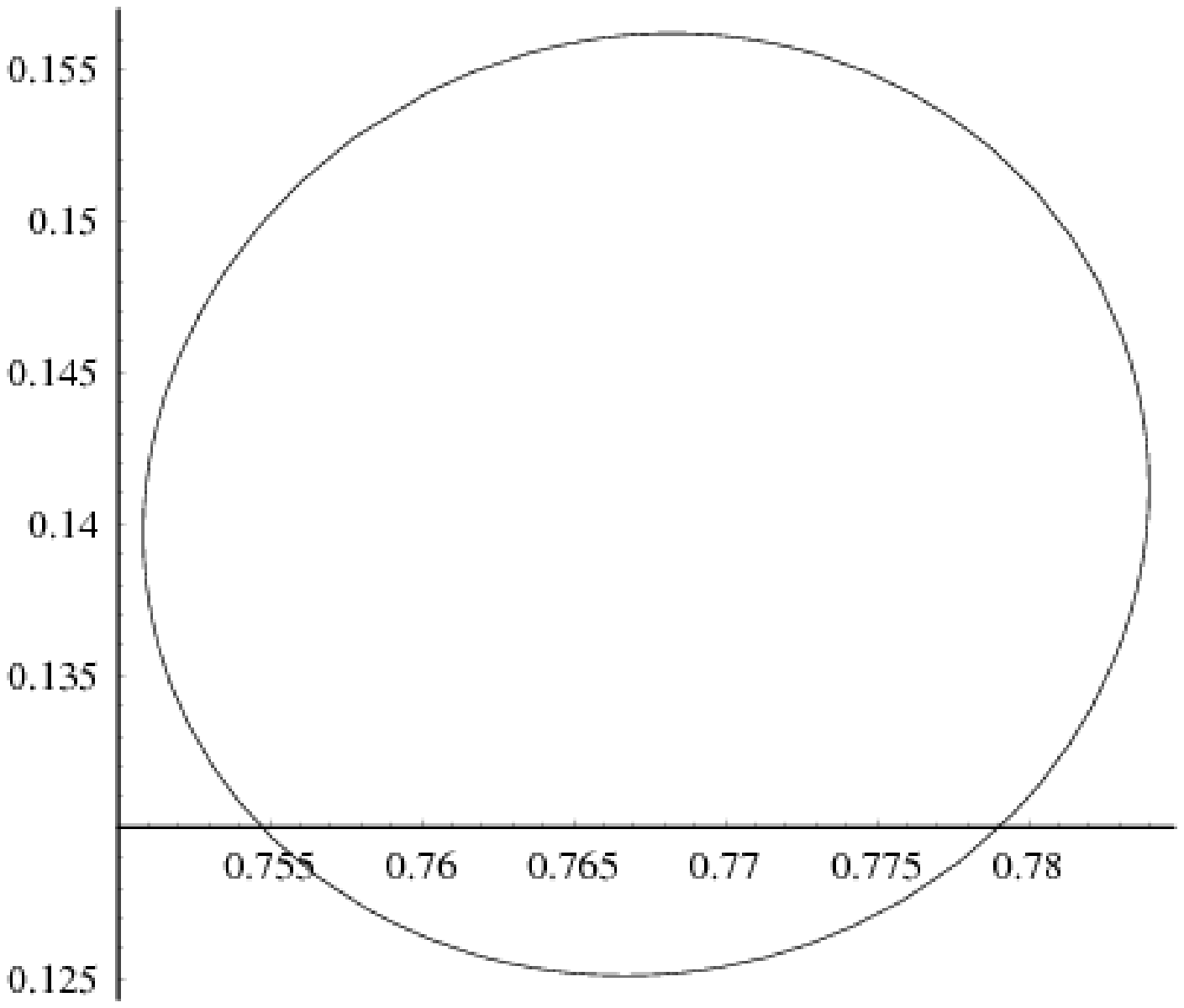}
\end{center}
\hspace{1.5cm}$\partial V_{\mathcal P}(z_0,\lambda)$ \hspace{4.5cm}  $\partial
V_{\mathcal G}(z_0,\lambda)$
\caption{´$z_0 =0.771264+0.151204i$ and $\beta=0.928608$}
\end{figure}

\begin{figure}[htp]
\begin{center}
\includegraphics[width=6.6cm]{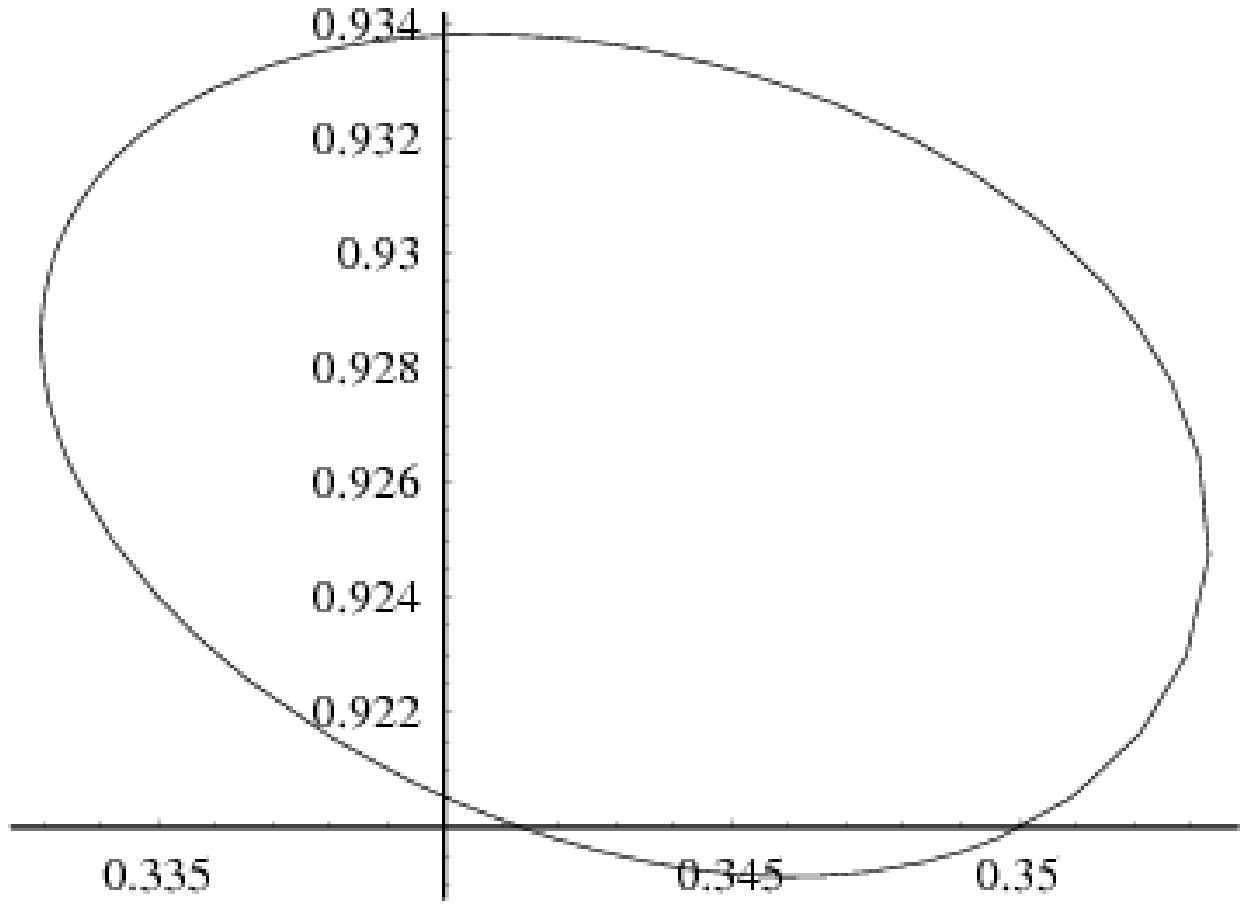}
\hspace{-0.5cm}
\includegraphics[width=6.6cm]{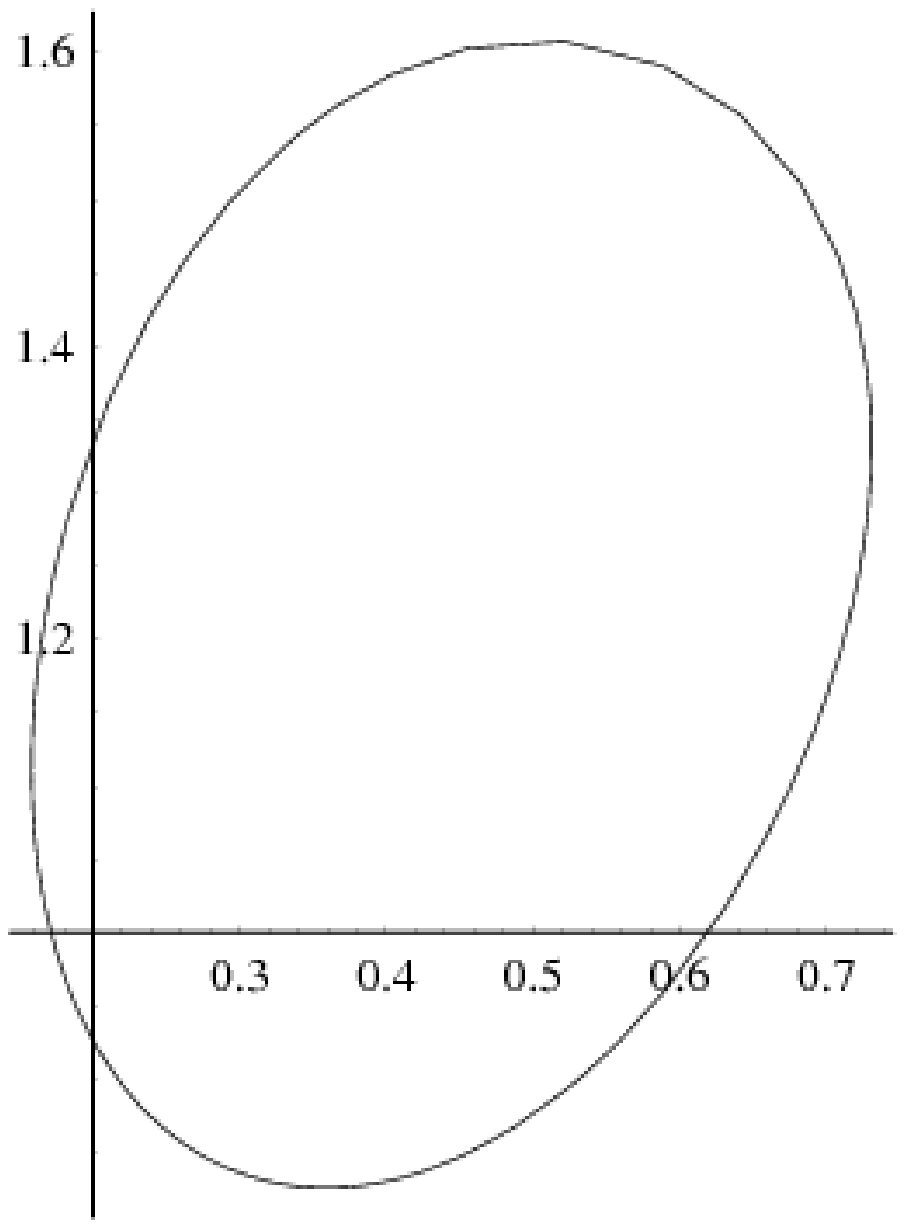}
\end{center}
\hspace{1.5cm}$\partial V_{\mathcal P}(z_0,\lambda)$ \hspace{4.5cm}  $\partial
V_{\mathcal G}(z_0,\lambda)$
\caption{$z_0 =0.335626+0.929093i$ and  $\beta =0.76622$}
\end{figure}

\vspace{1cm}

\subsection*{Acknowledgements}
This research work of the authors were supported by National Board for Higher Mathematics(DAE, India;
grant No. 48/2/2006/R\&D-II).

\empty\empty

\end{document}